\documentclass[reqno]{amsart}
\usepackage[utf8]{inputenc}
\usepackage{braket}
\usepackage{pmboxdraw}
\usepackage{amsmath}
\usepackage{amssymb}
\usepackage{xcolor}
\usepackage{graphicx}
\usepackage{float}
\usepackage[hidelinks]{hyperref}
\usepackage[shortlabels]{enumitem}

\calclayout

\usepackage{cite}
\usepackage{geometry}

\newtheorem{theorem}{Theorem}[section]

\newtheorem{lemma}[theorem]{Lemma}
\newtheorem{definition}[theorem]{Definition}
\newtheorem{proposition}[theorem]{Proposition}

\newtheorem{corollary}[theorem]{Corollary}

\newcommand{\N}{\mathbb{N}}

\newcommand{\R}{\mathbb{R}}

\newcommand{\ler}[1]{\left( #1 \right)}

\newcommand{\lers}[1]{\left\{ #1 \right\}}

\newcommand{\be}{\begin{equation}}
\newcommand{\ee}{\end{equation}}
\newcommand{\ba}{\begin{array}}
\newcommand{\ea}{\end{array}}

\newcommand{\cF}{\mathcal{F}}

\newcommand{\cC}{\mathcal{C}}

\newcommand{\cW}{\mathcal{W}}

\newcommand{\dd}{\mathrm{d}}

\newcommand{\Wp}{\mathcal{W}_p}
\newcommand{\cPX}{\mathcal{W}_p(X,d)}

\DeclareMathOperator*{\argmax}{arg\,max}
\DeclareMathOperator*{\argmin}{arg\,min}

\date{}



\begin{document}
\title{A generic $p$-Wasserstein space is isometrically rigid
}

\author[Zolt\'an M. Balogh]{Zolt\'an M. Balogh}
\address{Zolt\'an M. Balogh, Universit\"at Bern\\ Mathematisches Institut (MAI)\\ Sidlerstrasse 12\\ 3012 Bern\\ Schweiz}
\email{zoltan.balogh@unibe.ch}

\author[Eric Str\"oher]{Eric Str\"oher}
\address{Eric Str\"oher, Universit\"at Bern\\ Mathematisches Institut (MAI)\\ Sidlerstrasse 12\\ 3012 Bern\\ Schweiz}
\email{eric.stroeher@unibe.ch}

\author[Tam\'as Titkos]{Tam\'as Titkos}
\address{Tam\'as Titkos, Corvinus University of Budapest, Department of Mathematics \\
Fővám tér 13-15 \\ Budapest 1093 \\ Hungary\\ and HUN-REN Alfr\'ed R\'enyi Institute of Mathematics\\ Re\'altanoda u. 13-15.\\
Budapest 1053\\ Hungary}
\email{tamas.titkos@uni-corvinus.hu}

\author[D\'aniel Virosztek]{D\'aniel Virosztek}
\address{D\'aniel Virosztek, HUN-REN Alfr\'ed R\'enyi Institute of Mathematics\\ Re\'altanoda u. 13-15.\\Budapest H-1053\\ Hungary}
\email{virosztek.daniel@renyi.hu}

\thanks{Z. M. Balogh and E. Str\"oher are supported by the Swiss National Science Foundation, Grant Nr. {200020\_228012}.}

\thanks{T. Titkos is supported by the Hungarian National Research, Development and Innovation Office (NKFIH) under grant agreements no. K134944 and no. Excellence\_151232, and by the Momentum program of the Hungarian Academy of Sciences under grant agreement no. LP2021-15/2021.}

\thanks{D. Virosztek is supported by the Momentum program of the Hungarian Academy of Sciences under grant agreement no. LP2021-15/2021, by the Hungarian National Research, Development and Innovation Office (NKFIH) under grant agreement no. Excellence\_151232, and partially supported by the ERC Synergy Grant No. 810115.}

\begin{abstract}
We prove for any parameter $p \geq 1$ that the $p$-Wasserstein space $\cW_p(X,d_X)$ of a generic compact metric space $(X,d_X)$ is isometrically rigid, that is, its isometries are precisely the push-forwards of the isometries of $(X,d_X).$
We also study the question whether the isometry class of a $p$-Wasserstein space completely determines the isometry class of the underlying space. We prove for any $p\geq 1$ that for a generic compact metric space $(X,d_X),$ its $p$-Wasserstein space $\cW_p(X,d_X)$ is isometric to another $p$-Wasserstein space $\cW_p(Y,d_Y)$ if and only if $(Y,d_Y)$ is isometric to $(X,d_X).$ En route we prove that a stronger property, which we call the Dirac preservation property, is also generic.
\end{abstract}

\maketitle

\section{Introduction}

\subsection{Motivation and main results}

The geometry of Wasserstein spaces and, in particular, the structure of their isometries, has been intensively studied in recent years. The isometry group of a $p$-Wasserstein space $\Wp(X, d)$ (see the precise definition is Subsection \ref{susec:basic-not-not}) depends on the underlying metric space $(X, d)$ and also on the value of the parameter $p \geq 1.$
It is known that the push-forward operation embeds the isometry group of $(X,d)$ into that of $\Wp(X,d)$. More precisely, if $\varphi:X\to X$ is an isometry of $X$, then $\varphi_{\#}$ is an isometry of $\Wp(X,d)$. Isometries of the form $\varphi_{\#}$ are thus called trivial. A fundamental question that has been extensively studied in recent years, see e.g. \cite{BKTV, BTV, BTV2,BSTV, BK2,BK,GTV1,GTV2,K,S-R,S-R2}, is whether isometries of other types exist. 
If the isometry group of $\Wp (X, d)$ contains only trivial isometries, we say that the Wasserstein space is \emph{isometrically rigid}; otherwise, we call the space non-rigid or flexible. 

While there are remarkable examples of flexible Wasserstein spaces, \cite{K, S-R2, GTV1} it seems that the vast majority of Wasserstein spaces are isometrically rigid, that is, its isometry group coincides with that of the underlying space: \cite{BK, BK2, BSTV, BKTV, BTV, BTV2}. Moreover, in \cite{BSV} it was shown, that adding just one point to the known examples metric spaces for which the Wasserstein space is non-rigid the rigidity is restored.  Furthermore, in \cite{BSV} the question has been raised whether the rigidity of the Wasserstein space is a generic property in the topology of the Gromov-Hausdorff convergence of compact metric spaces. For the general background on Gromov-Hausdorff convergence, we refer the reader to \cite{Bate,BBI,Herron}.

The first main goal of the present work is to answer the question from \cite{BSV} in the affirmative. This result, stated as Theorem \ref{thm:rigidity-is-generic} below, formalizes the experience that most Wasserstein spaces are isometrically rigid as a precise mathematical statement.

\begin{theorem} \label{thm:rigidity-is-generic}
For any fixed parameter $p \geq 1,$ the set of compact metric spaces $(X, d_X)$ with isometrically rigid $p$-Wasserstein  spaces $\Wp(X, d_X)$ form a residual set in the  Gromov-Hausdorff topology.   
\end{theorem} 

We note that the statement of Theorem \ref{thm:rigidity-is-generic} follows for $p=1$ from the result of Rouyer \cite[Theorem 3]{Rouyer} saying that in a generic compact metric space three distinct points are never collinear, and from \cite[Theorem 4.6]{GTV2} saying that the $1$-Wasserstein space over a metric space with no collinear triple of points is isometrically rigid.
\par 
The proof of this theorem has two distinct parts. In the first part we prove that for all $p\geq 1$ the set of compact  metric spaces $X$ having rigid Wasserstein space $\mathcal{W}_p(X,d_X)$ is dense in the Gromov-Hausdorff topology of all compact metric spaces. In the second part we prove that the same set: i.e. the set of all compact metric spaces $X$ for which the Wasserstein space $\mathcal{W}_p(X,d_X)$ is rigid is a countable intersection of open sets, that is, a $G_\delta$ set. Combining these statements we obtain that compact metric spaces $X$ having rigid Wasserstein space $\mathcal{W}_p(X,d_X)$ form a dense $G_{\delta}$ set. Since dense $G_{\delta}$ sets in any complete metric space are residual, the result follows. 
\par
We also study the question whether the isometry class of a $p$-Wasserstein space completely determines the isometry class of the underlying space. That is, does $\cW_p(X, d_X) \simeq \cW_p(Y,d_Y)$ imply that $(X, d_X) \simeq (Y,d_Y)$? (Here and throughout this paper, the symbol $\simeq$ stands for the isometric relation between metric spaces.) Bertrand and Kloeckner answered this question for $p=2$ in the special case of locally compact and geodesically complete spaces affirmatively \cite[Theorem 1.2]{BK}, and Lenze recently provided an affirmative answer also for complete connected Riemannian manifolds \cite[Theorem C]{Lenze} in the $p=2$ case. However, in full generality, this question remained open. We have the following positive result in this direction: for any $p \geq 1,$ a generic compact metric space $(X,d_X)$ satisfies the recovery property of order $p$ (abbreviated as RECp, see Definition \ref{def:RECp}), that is, its $p$-Wasserstein space $\cW_p(X,d_X)$ is isometric to another $p$-Wasserstein space $\cW_p(Y,d_Y)$ if and only if $(Y,d_Y)$ is isometric to $(X,d_X).$ 
We actually prove a stronger statement: we show that for all $p \geq 1,$ a generic compact space satisfies the so-called Dirac preservation property of order $p$ (abbreviated as DPPp, see Definition \ref{def:DPPp}). We also demonstrate that the Dirac preservation property of order $p$ is strictly stronger than the recovery property of order $p$ for all $p \geq 1.$ The precise statement, which is the second main result of this paper, is formalized below. 

\begin{theorem} \label{thm:recovery-is-generic}
For any $p\geq 1$ the set of metric spaces $X \in \mathfrak{M}$ satisfying the Dirac preservation property DPPp is a residual set in $\mathfrak{M}.$ Consequently, a generic compact metric space satisfies the recovery property RECp.
\end{theorem}

%==============================================
\subsection{Basic notions, notation} \label{susec:basic-not-not}

For a fixed $p\geq1$ and a complete separable metric space $(X, d)$, the $p$-Wasserstein space $\Wp(X, d)$ is the set of Borel probability measures with finite moment of order $p,$ equipped with the so-called Wasserstein distance $d_{\Wp}$. More precisely, we say that a Borel probability measure $\mu$ is in the $p$-Wasserstein space if 
$$
\int_X d^p(x, x_0)~\mathrm{d}\mu(x)<\infty
$$
for some (and hence any) $x_0 \in X$, and the Wasserstein distance is given by 
$$
d_{\mathcal{W}_p}(\mu, \nu)=\inf_{\pi \in \Pi (\mu, \nu)}\left\lbrace \left( \iint_{X \times X} d_X^p(x, y)~\mathrm{d}\pi(x,y)\right)^{\frac{1}{p}} \right\rbrace, 
$$
where $\Pi(\mu, \nu)$ is the set of couplings between $\mu$ and $\nu$, i.e. the set of probability measures $\pi$ on $X \times X$ with $p_{1 \#}\pi=\mu$ and $p_{2 \#}\pi=\nu$, for the projection operators $p_1(x, y)=x$ and $p_2(x, y)=y$; here $p_{\#}$ denotes the push-forward operation, which we define below. If $\pi$ is a coupling between $\mu$ and $\nu$, we denote its Kantorovich cost by $K(\pi)$,
$$K(\pi) = \int_{X\times X} d^p(x,y) ~\mathrm{d}\pi(x,y).$$
We say that the coupling $\pi$ is optimal if 
$K(\pi) = d^p_{\mathcal{W}_p}(\mu, \nu)$.

The Wasserstein space $\Wp(X, d)$ inherits a number of properties from its base space $(X, d)$, see \cite{ Villani} for an overview. For example, if $(X, d)$ is a Polish metric space, then so is $\Wp(X, d).$ 
Another property of any Wasserstein space is that it contains an isometric copy of its base space. Indeed, the map $x \mapsto \delta_x$ is an isometric embedding from the space $(X, d)$ to $\Wp(X, d),$ where $\delta_x$ is the Dirac mass located at $x \in X$. Furthermore, the convex hull of Dirac masses (i.e. the set of measures with finite support) is dense in $\Wp(X,d)$ (see e.g. \cite[Theorem 6.18]{Villani}). The push-forward operation lifts a measurable map $T: (X,d) \to (X,d)$ to a transformation $T_\#$ of the $p$-Wasserstein space $\Wp(X,d)$ which is defined by 
$$
T_{\#}(\mu)(A) = \mu (T^{-1}(A))\quad\mbox{for all measurable subset } A \text{ of } X.
$$
Throughout this work, by \emph{isometry} we mean a distance-preserving and bijective map on a metric space.
\par
Let us recall now the Gromov-Hausdorff topology (see Chapter 7 of \cite{BBI}) on spaces of metric spaces. First, let us consider the case when two compact sets $A, B \subset Z$ where $(Z, d)$ is a given metric space. In this case the Hausdorff distance $d_H(A, B)$ is defined as 
$$d_H(A,B) = \inf \{ r> 0: A \subseteq U_r(B) \ \text{and} \ B\subseteq U_r(A) \}, $$
where $U_r(A)= \{ y \in Z: \exists \  x \in A \ \text{with} \ d(x,y) <r \} $ is the $r-$ neighborhood of $A$ and similarly for $U_r(B)$. It turns out that $d_H$ is a metric on the space $\mathfrak{M}(Z)$ of all compact subsets of $Z,$ which makes $\mathfrak{M}(Z)$ compact whenever $Z$ is compact. 

Using this, we can define the Gromov-Hausdorff distance
$d_{GH}(X,Y)$ between compact metric spaces $X, Y$ as the inf of all $r>0 $ with the property that there exists a metric space $Z$ and compact subsets of $X', Y' \subseteq Z$ such that $X$ is isometric to $X'$, $Y$ is isometric to $Y'$ and $d_H(X', Y') < r$. Note, that this distance will not be a metric on the spaces of all compact metric spaces $\mathfrak{M}$ since the distance between to distinct but isometric spaces is $0$. However it will be a metric on the isometry classes of compact metric spaces. To ease notation we shall replace isometry classes by considering a representative in every class. In this sense we say that a sequence of compact metric spaces $(X_n)_{n\in\mathbb{N}}$ converges in the Gromov-Hausdorff distance to an element $X \in \mathfrak{M}$ if $d_{GH}(X_n, X) \to 0$ when $n\to \infty$. 
Recall that a subset of a topological space is said to be rare (in other words: nowhere dense) if the interior of its closure is empty. It is said to be meager, or of first category, if it is a countable union of rare subsets. A space is called a Baire space if each meager subset of it has empty interior. The Baire-Hausdorff theorem states that any complete metric space is a Baire space. The complement of a meager subset of a Baire space is said to be residual. At last, given a Baire space, we say that a generic element enjoys a property if the set of elements which satisfy this property is residual.

\section{Spaces with rigid $p$-Wasserstein spaces form a dense subset of all compact metric spaces}
As described in the previous section, we  consider compact metric spaces $(X, d)$ endowed with the Gromov-Hausdorff topology. Following \cite{Rouyer} we call a metric space  $(X, d)$  totally anisometric if two distinct pairs of points have distinct distances. Let us recall that by Corollary 4 from \cite{Rouyer} finite metric spaces that are totally anisometric are dense in the space $\mathfrak{M}$ of all compact metric spaces $(X, d)$ with the Gromov-Hausdorff topology. We consider a subset of finite, totally anisometric spaces defined as follows:

\begin{definition} \label{def:NDp}
Let $p\geq 1$ and a finite metric space $ (X, d)$ where $X = \{ x_1, \ldots, x_N \}$. We say that $(X, d)$  satisfies nondegenerate property (NDp) if the $N\times N$ matrix defined by the distances of pairs of points raised to the power $p$ has non-vanishing determinant:
\begin{equation} \label{nondeg}
\det [d^p(x_i, x_j)] \neq 0.
    \end{equation}
\end{definition}
We note that the set of finite, totally anisometric spaces satisfying (NDp) are dense in $\mathfrak{M}$. The density of finite, totally anisometric spaces in $\mathcal{M}$ is clear, and hence the only thing we have to show here is that every such space can be perturbed to satisfy the NDp as well. Fix $p\geq 1$, and let $(X=\{x_1,\ldots,x_N\},d)$ be a finite, totally anisometric space. For every $t>0,$ consider the modified metric
\begin{align}
d_t(x_i,x_j)=
\begin{cases}
\ler{d^p(x_i,x_j)+t}^{1/p},&i\ne j,\\
0,&i=j,
\end{cases} \nonumber
\end{align}
which is the $l_p$ sum of $d$ and the rescaled discrete metric. The function
$ r\mapsto(r^p+t)^{1/p}$ is strictly increasing, so distinct distances remain distinct. Write
$$
D=\bigl[d^p(x_i,x_j)\bigr]_{i,j=1}^N, \, D_t=\bigl[d_t^p(x_i,x_j)\bigr]_{i,j=1}^N,
$$
and observe that $D_t=D+t(J-I)$ where $J$ is the all-ones matrix and $I$ is the identity matrix.
Consequently,
$$
P(t):=\det\bigl(D+t(J-I)\bigr)
$$
is a polynomial in $t.$ Its leading coefficient is
$$
\det(J-I)=(N-1)(-1)^{N-1}\ne0,
$$
because $J-I$ has eigenvalues $N-1$ and $-1$, the latter with multiplicity $N-1.$
Thus $P$ is not the zero polynomial and has only finitely many roots. In particular, there are arbitrarily small positive $t$ for which
$$
\det D_t\neq 0.
$$
On the other hand, it is clear that $(X,d_t)$ converges to $(X,d)$ in the Gromov-Hausdorff distance as $t \to 0.$
\par
The main result of this section is the following:
\begin{proposition} \label{finite-rigid}
  Let $p \geq 1,$ and let $(X, d)$ be a finite, totally anisometric metric space satisfying (NDp).  Then that the $p$-Wasserstein  spaces $\Wp(X, d)$ is isometrically rigid.    
\end{proposition}

Let us note that the isometry group, $\mathrm{Isom}(X)$ of finite,  totally anisometric spaces $X$ reduces to the identity map $\mathrm{Id}_X$ if $|X| \geq 3.$ Therefore the the above statement implies that the isometry group $\mathrm{Isom}(\Wp(X, d))$ for such spaces reduces to the identity map $\mathrm{Id}_{\Wp(X,d)}$ as well. Since we shall use this observation later, we formulate it as the following: 

\begin{corollary} \label{rigid-dense}
 For any parameter $p \geq 1,$ the subset $\mathcal{R}_p\subset \mathcal{M}$ of compact metric spaces $(X, d)$ with the property that $\mathrm{Isom(}\Wp(X, d))= \{\mathrm{Id}_{\Wp(X,d)} \}$ is dense in the spaces of compact metric spaces $\mathfrak{M}$ with respect to the Gromov-Hausdorff topology.       
\end{corollary}

The proof of Proposition \ref{finite-rigid} is based on the following:

\begin{proposition} \label{prop-dirac-preservation}
    For any parameter $p \geq 1,$ and any finite, totally anisometric metric space $(X, d_X)$ containing at least three points, it follows that if  $\Phi: \Wp(X, d)\to \Wp(X, d)$ is an isometry, then it preserves the Dirac masses i.e.
    $\Phi(\delta_x) = \delta_x $ for all $x \in X$.    
\end{proposition}

In the following we indicate the short proof of Proposition \ref{finite-rigid} based on Proposition \ref{prop-dirac-preservation}.

\begin{proof}[Proof of Proposition \ref{finite-rigid}]
Let $p\geq 1$ and $X = \{ x_1, \ldots, x_N \}$ be an anisometric,  finite metric space satisfying (NDp) and let $\Phi: \Wp(X, d)\to \Wp(X, d)$ be an isometry. From Proposition \ref{prop-dirac-preservation} we obtain that $\Phi(\delta_{x_i}) = \delta_{x_i}$ for all $i= 1, \ldots, N$. 

Let us consider an arbitrary element 
$\mu \in \Wp(X, d)$ written as $\mu = \sum_{i=1} ^N \lambda_i\delta_{x_i} $ with some coefficients $0\leq \lambda_i \leq 1$ with $\sum_i \lambda_i = 1$. Denoting by $\mu' = \Phi(\mu)$ we can write $\mu' = \sum_{i=1} ^N \lambda'_i\delta_{x_i} $ with some coefficients $0\leq \lambda'_i \leq 1$ with $\sum_i \lambda'_i = 1$. 

Our goal is to show that $\lambda_i = \lambda'_i$ for all $i= 1, \ldots, N$. To do that we use Proposition \ref{prop-dirac-preservation}. From the fact that $\Phi: \Wp(X, d)\to \Wp(X, d)$ is an isometry, it follows
$$d^p_{\Wp}(\mu, \delta_{x_j}) = d^p_{\Wp}(\Phi(\mu), \Phi(\delta_{x_j})) =  d^p_{\Wp}(\mu', \delta_{x_j}) , \ j = 1, \ldots, N . $$

Observing that for all $j = 1, \ldots, N $ between $\mu$ and $\delta_{x_j}$ the only coupling is the trivial one, we can write: 
  $$ d^p_{\Wp}(\mu, \delta_{x_j})  = \sum_{i=1}^N\lambda_i d^p(x_i, x_j)  \ \text{and} \  d^p_{\Wp}(\mu', \delta_{x_j})  = \sum_{i=1}^N\lambda'_i d^p(x_i, x_j) \ ,\ j= 1, \ldots, N  $$
  which shows that the vector $(\lambda_1- \lambda_1', \ldots, \lambda_N-\lambda_N')$ is a solution to the system of linear equations:
  $$ \sum_{i=1} ^N (\lambda_i -\lambda'_i)  d^p(x_i, x_j) =0 \ , \ j= 1, \ldots, N .  $$

By condition (NDp) the coefficient matrix of the above system is invertible and the only solution of the system is the trivial one, which gives $\lambda_i = \lambda'_i$ for all $i= 1, \ldots, N$, finishing the proof. 
\end{proof}

\subsection{Proof of Proposition \ref{prop-dirac-preservation}}

In this subsection we shall give the proof of Proposition \ref{prop-dirac-preservation} that will be based on a metric characterization of Dirac measures in $\Wp(X, d)$ formulated in Lemma \ref{lemma:metric-ratio-singleton} below. 

Before stating this Lemma we shall start with two definitions. The first one defines a relation on $X.$ Let $\mathfrak{F} \subset X \times X$ be defined in the following way: 
\begin{align} \label{eq:frak-F-def}
    \mathfrak{F}:=\left\{(v,v') \in X \times X \, \middle| \, d^p(v,v') < \sum_{l=0}^{L-1} d^p(x_l, x_{l+1})  \text{ for any finite chain }\right. \nonumber \\ 
     \left.  (v=x_0,x_1, x_2, \dots, x_{L-1},x_L=v') \neq (v,v') \text{ of distinct points in } X\right\}.
\end{align}
Note first that $\mathfrak{F}$ is a symmetric relation on $X \times X$ and hence we can associate an undirected graph $G=(X,\mathfrak{F})$ to it with vertex set $X$, and declaring that an edge between $(x,x')$ exists if and only if $(x,x')\in \mathcal{F}$.

\begin{lemma} \label{lemma:X-F-graph-connected} 
    Let $p \geq 1$, and let $(X,d)$ be a finite, totally anisometric metric space with at least three points. Then the relation $\mathfrak{F} \subseteq X \times X$ defined by \eqref{eq:frak-F-def} gives rise to a connected graph on $X,$ that is, the graph $G=(X,\mathfrak{F})$ is connected. 
\end{lemma}

\begin{proof}

    To prove the statement of the claim we assume by contradiction that our graph has more than one connected component, and let $(z, z')$ be the shortest edge between different connected components. In mathematical terms: 
    \begin{align}
        (z, z'):=\argmin \lers{d^p(x,x') \, \middle| \, x \text{ and } x' \text{ are in different connected components of}~ (X,\mathfrak{F})} \nonumber
    \end{align}
    But then $(z, z') \in \mathfrak{F}$. Indeed, any path $(z=x_0, x_1, \dots, x_M=z')$ connecting $z$ with $z'$ necessarily contains a step between two different connected components. Let $(x_m,x_{m+1})$ be this step. By the minimality of $d^p(z,z')$ we have $d^p(z,z') < d^p(x_m,x_{m+1}),$ and hence $d^p(z,z') < \sum_{j=0}^{M-1}d^p(x_j, x_{j+1}),$ which means by the definition \eqref{eq:frak-F-def} that $(z,z') \in \mathfrak{F}.$ This is a contradiction, proving the claim.  
\end{proof}

The second definition reads as follows:
If $ \Lambda:=(\lambda_i)_{i=1}^N \in [0,1]^N$ for $N \geq 1$ is a set of strictly increasing weights, we consider the multiple ratio $M_{\Lambda}(\mu, \nu)\subset \cPX^N$ defined by
\begin{align}
M_{\Lambda}(\mu, \nu)=\lbrace\{\xi_i\}_{i=1}^{N} \in \cPX^N | 
 d_{\mathcal{W}_p}(\xi_i, \xi_{i+1})\leq (\lambda_{i+1}-\lambda_i)^{1/p}d_{\mathcal{W}_p}(\mu, \nu), 0 \leq i \leq N+1 \rbrace. \label{eq:multiple-ratio-set-def}
\end{align} using the notation that $\lambda_0=0,$ $\lambda_{N+1}=1$ and $\xi_0=\mu, $ $\xi_{N+1}=\nu$. 
When $\Lambda=\{\lambda\}$ is a single element, we write $M_{\lambda}(\mu, \nu):=M_{\Lambda}(\mu, \nu)$. 

Let us notice that for any two measures $\mu, \nu \in \cPX$ the set $M_{\Lambda}$ is invariant under isometry, that is, if $\{\xi_i\}_{i=1}^N \in M_{\Lambda}(\mu, \nu)$, then $\{\Phi(\xi_i)\}_{i=1}^N \in M_{\Lambda}(\Phi(\mu), \Phi(\nu))$. This property follows immediately from the distance preserving property of isometries and the definition of $M_{\Lambda}(\mu, \nu)$. In particular, since isometries are bijective, $M_{\Lambda}(\mu, \nu)$ is a singleton if and only if $M_{\Lambda}(\Phi(\mu), \Phi(\nu))$ is a singleton.

The following simple claim asserts the existence of at least one element in $M_{\Lambda}(\mu, \nu)$, as the appropriate sequence of standard convex combinations of $\mu$ and $\nu$ always belongs to it.

\begin{lemma} \label{lemma:existence of Mlambda}
Let $p \geq 1$, and let $(X,d)$ be a finite, totally anisometric metric space with at least three points. For any $N \geq 1 $, any two measures $\mu, \nu \in \cPX$ and any $\Lambda =(\lambda_i)_{i=1}^N\in [0,1]^N$, we have 
$$
\{\xi_i\}_{i=1}^N = \{(1-\lambda_i) \mu + \lambda_i \nu\}_{i=1}^N \in M_{\Lambda}(\mu, \nu).
$$
\end{lemma}

\begin{proof}
Consider $\pi_0$ the optimal transport plan between $\mu$ and $\nu$, and set 
$$
\pi_{i+1}= (1-\lambda_{i+1})(\mathrm{Id}_X \times \mathrm{Id}_X)_{\#}\mu + \lambda_{i}(\mathrm{Id}_X \times \mathrm{Id}_X)_{\#}\nu +  (\lambda_{i+1}-\lambda_i) \pi_0,
$$ with $\lambda_0=0$. 
Then $\pi_{i+1}$ is a transport plan between $\xi_{i} = (1-\lambda_i) \mu + \lambda_i \nu$ and $\xi_{i+1}$  for any $0 \leq i \leq N$ (with the convention that $\xi_0=\mu, \xi_{N+1}=\nu, \lambda_0=0, \lambda_{N+1}=1$).
Indeed, it is easy to see that the left marginal of $\pi_{i+1}$ is
$$(1-\lambda_{i})\mu + \lambda_i\nu,$$ while the right marginal is 
$$
(1-\lambda_{i+1})\mu + \lambda_{i+1}\nu.
$$
The Kantorovich cost is given by $K(\pi_{i+1})=(\lambda_{i+1}-\lambda_i) K(\pi_0)$, which implies that $$d_{\mathcal{W}_p}(\xi_i, \xi_{i+1} ) \leq (\lambda_{i+1}-\lambda_i)^{1/p}d_{\mathcal{W}_p}(\mu, \nu),$$ giving that $\{\xi_i\}_{i=1}^N \in M_{\Lambda}(\mu, \nu)$. 
\end{proof} 

In the following, whenever we consider a set $\Lambda \in [0,1]^N$, we will assume that it is strictly increasing. 
We now prove the following lemma, characterizing those pairs of measures $(\mu, \nu)$ with the property that which $M_{\Lambda}(\mu, \nu)$ contains the unique sequence of measures given by their convex combination. This result will be the important technical tool in proving Proposition \ref{prop-dirac-preservation}:

\begin{lemma}  \label{lemma:metric-ratio-singleton}
Let $p \geq 1$, and let $(X,d)$ be a finite, totally anisometric metric space with at least three points, and let $\mu, \nu\in\cPX$. Then the following are equivalent
\begin{itemize}
\item[(a)]  There exists $(v, v') \in \mathfrak{F}$, $\eta \in \cPX$ and $c \in [0,1]$ such that $$\mu=(1-c) \eta + c \delta_v\qquad\nu=(1-c)\eta + c \delta_{v'}.$$
\item[(b)] For any $\Lambda=(\lambda_i)_{i=1}^N\subset [0,1]$ set of strictly increasing weights, the multiple ratio set $M_{\Lambda}(\mu, \nu)$ is a unique set of $N$ measures given by $M_{\Lambda}(\mu, \nu)=\{(1-\lambda_i)\mu + \lambda_i \nu\}_{i=1}^N$.
\end{itemize}
\end{lemma}

In what follows we shall give the proof of Proposition \ref{prop-dirac-preservation} based on Lemma \ref{lemma:metric-ratio-singleton} and Lemma \ref{lemma:X-F-graph-connected}: 

\begin{proof}[Proof of Proposition \ref{prop-dirac-preservation}] 
Consider an isometry $\Phi:\cPX \to \cPX$ and assume that there exists an $x_0 \in X$ such that $\Phi(\delta_{x_0})=\delta_{x_0}$. In the first step we show that for any $x\in X$ such that $(x_0, x) \in \mathfrak{F}$, we have that $\Phi(\delta_{x})=\delta_x$. The proof of this fact will be based on Lemma \ref{lemma:metric-ratio-singleton}.

Now consider the set $F \subset X$ where $x \in F$ if $(x_0, x) \in \mathfrak{F}$. Since $X$ is finite, so is $F$, and we can write $F=\{x_i\}_{i=1}^K$, for some $K\geq 1 $  with the considered order that $d(x_0, x_i)>d(x_0, x_j)$ whenever $i<j$ (that is, the values are ordered in the inverse order of their distance from $x_0$). Note that here we use the assumption that $X$ is totally anisometric. 

We apply Lemma \ref{lemma:metric-ratio-singleton} in the following way: pick $c:=1$, and write $\delta_{x_0} =(1-1)\eta + \delta_{x_0}$ and $\delta_{x}=(1-1)\eta + \delta_{x}$, we have that if $(x_0, x) \in \mathfrak{F}$, then by the statement of Lemma \ref{lemma:metric-ratio-singleton} we obtain that the set $M_{\Lambda}(\delta_{x_0}, \delta_x)$ is a singleton for any $\Lambda \in [0,1]^N$. Thus $M_{\Lambda}(\Phi(\delta_{x_0}), \Phi(\delta_x))$ will also be a singleton, and using Lemma \ref{lemma:metric-ratio-singleton} again there exists $\eta\in \cPX$, $c \in [0,1]$ such that $\Phi(\delta_{x_0})=(1-c)\eta + c \delta_{v}$ and $\Phi(\delta_x)=(1-c)\eta + c \delta_{v'}$ with $(v, v') \in \mathfrak{F}$. Since we assumed that $\Phi(\delta_{x_0})=\delta_{x_0}$, we have that either $c=1$ or $\eta=\delta_{x_0}$ and in any of these two cases we must have in addition that $v=x_0$. Thus we have $\Phi(\delta_{x})=(1-c)\delta_{x_0} + c \delta_{x'}$ with $(x_0, x') \in \mathfrak{F}$. 

We first show that $\Phi(\delta_{x_1})=\delta_{x_1}$. To see this, note, that by the previous consideration, we have that $\Phi(\delta_{x_1})=(1-c)\delta_{x_0} + c \delta_{x_i}$ for some $1 \leq i \leq K$. But we also have that 
$$d_{\mathcal{W}_p}(\Phi(\delta_{x_0}), \Phi(\delta_{x_1}))=d_{\mathcal{W}_p}(\delta_{x_0}, (1-c)\delta_{x_0} + c \delta_{x_i})= c^{1/p}d(x_0, x_i).$$ Thus, since $\Phi$ preserves distances, 
$$
d(x_0, x_1)= d_{\mathcal{W}_p}(\delta_{x_0}, \delta_{x_1})=d_{\mathcal{W}_p}(\Phi(\delta_{x_0}), \Phi(\delta_{x_1}))=c^{1/p}d(x_0, x_i).
$$
As by the assumption $d(x_0, x_1)>d(x_0, x_i)$ whenever $i>1$, and that $c\leq 1$ we conclude that in fact $c=1$ and $x_i=x_1$, and thus $\Phi(\delta_{x_1})=(1-1)\eta + \delta_{x_1}=\delta_{x_1}$. 

As a consequence of this, we also immediately get that for any $c \in [0,1]$, the interpolation measure $\xi=(1-c)\delta_{x_0}+ c \delta_{x_1}$ is also preserved under isometry, $\Phi(\xi)=\xi$. Indeed, for $\Lambda=\{c\}$, $M_{\Lambda}(\delta_{x_0}, \delta_{x_1})=\{\xi\}$, and thus $\{\Phi(\xi)\}=M_{\Lambda}(\Phi(\delta_{x_0}), \Phi(\delta_{x_1}))=M_{\Lambda}(\delta_{x_0}, \delta_{x_1})=\{\xi\}$. 
 
We now finish the proof by induction. To do so,  assume that for any $i < k$, $\Phi(\delta_{x_i})=\delta_{x_i}$ and for any $c \in [0,1]$, $\Phi((1-c)\delta_{x_0} + c \delta_{x_i})= (1-c)\delta_{x_0} + c \delta_{x_i}$. As in the previous parts, we have that $\Phi(\delta_{x_k})=(1-c)\delta_{x_0} + c \delta_{x_j}$. Since 
$$
d(x_0, x_k)=d_{\mathcal{W}_p}(\Phi(\delta_{x_0}),\Phi( \delta_{x_k}))=c^{\frac{1}{p}} d(x_0, x_j),
$$ we have that $j \leq k$, with equality if and only if $c=1$. If $j<k$, then by the induction assumption $$\Phi^{-1}((1-c)\delta_{x_0} + c \delta_{x_j})=(1-c)\delta_{x_0} + c \delta_{x_j} \neq \delta_{x_k},$$ and $\Phi(\delta_{x_k})=(1-1)\delta_{x_0} + \delta_{x_k}=\delta_{x_k}$. Furthermore, we have by the same method as above that $\Phi((1-c)\delta_{x_0} + c\delta_{x_k})=(1-c)\delta_{x_0} + c\delta_{x_k}$. Using this induction step $K-1$ times we have thus proven that $\Phi(\delta_{x_i})=\delta_{x_i}$ for any $x_i$ such that $(x_i, x_0)\in \mathfrak{F}$. 

 Now we can apply Lemma \ref{lemma:X-F-graph-connected}: by the fact that the graph induced by the relation $\mathfrak{F}$ is connected it shows that if we can find one single point $x_0 \in X $ such that $\Phi(\delta_{x_0})= \delta_{x_0}$ then by the previous step we can propagate the preservation of Dirac masses according to the relation to every point in $X$. This would imply that $\Phi(\delta_x) = \delta_x$ for any $x\in X$. 

 To finish the proof we must show that we can start somewhere with the preservation of the Dirac masses: i.e. there exists a point $x_0 \in X$ such that $\Phi(\delta_{x_0})= \delta_{x_0}$. 

To do so, we use the anisometric property of $X$ to find unique pair of points realizing the diameter in $(X,d)$ denoted by $\{x_0, x'_0\}.$ Then the unique pair of point realizing the diameter in $\cW_p(X,d)$ is $\{\delta_{x_0},\delta_{x'_0}\}.$ Consequently, $\{\Phi(\delta_{x_0}),\Phi(\delta_{x'_0})\}$ is also the unique pair of point realizing the diameter in $\cW_p(X,d),$ which implies that  $\{ \Phi(\delta_{x_0}), \Phi(\delta_{x'_0})\} = \{ \delta_{x_0}, \delta_{x'_0} \} $. If  $\Phi(\delta_{x_0}) = \delta_{x_0}$ then we are done. We show below that the other option $\Phi(\delta_{x_0}) = \delta_{x'_0}$ will lead to a contradiction.

To do so, we let $\sim$ denote the binary relation on $\cW_p(X,d)$ characterized by Lemma \ref{lemma:metric-ratio-singleton}, meaning that $\mu \sim \nu$ iff $\mu$ and $\nu$ are represented according to the condition $a)$ from 
Lemma \ref{lemma:metric-ratio-singleton}. Taking $\mu = \delta_{x_0}= (1-c)\delta_{x_0} + c\delta_{x_0} $ for some $c\in [0,1]$ we obtain that 

\begin{align} \label{eq:M-char}
    \mathcal{M}=\lers{\nu \in \cPX: \nu \sim \delta_{x_0} }=\lers{(1-c)\delta_{x_0} + c \delta_x \, \middle| \, (x_0 ,x) \in \cF, \ c \in [0,1]},
\end{align}
and similarly, 
\begin{align} \label{eq:M-prime-char}
    \mathcal{M'}=\lers{\nu' \in \cPX: \nu' \sim \delta_{x'_0} }=\lers{(1-c)\delta_{x'_0} + c \delta_x' \, \middle| \, (x'_0 ,x') \in \cF, \ c \in [0,1]}.
\end{align}
The relation $\sim$ is characterized in point (b) of Lemma \ref{lemma:metric-ratio-singleton} purely in terms of the $p$-Wasserstein metric, and hence $\mathcal{M'}=\{ \Phi(\nu) \, | \, \nu \in \mathcal{M}\}.$ Clearly, $d_{\Wp}(\delta_{x_0}, (1-c)\delta_{x_0} + c \delta_x)=c^{1/p} d(x_0,x),$ and hence $\max \{d_{\Wp}(\delta_{x_0}, \nu) \, | \, \nu \in \mathcal{M} \}=\max \{d(x_0,x) \,|\, (x_0,x) \in \cF\}$. In the same way we obtain that $\max \{d_{\Wp}(\delta_{x'_0}, \mu') \, | \, \mu' \in \mathcal{M'} \}=\max \{d_X(x'_0,x') \,|\, (x'_0,x') \in \cF\},$ and that these two quantities must be equal to each other by the isometric property of  $\Phi$  and the fact that $\Phi(\delta_{x_0})= \delta_{x'_0} $. By our assumption the space $X$ is anisometric. This implies  that 
$(x_0, x'_0) \in \mathfrak{F}$. 

Considering $\mu=\delta_{x_0} $ and $\nu= \delta_{x'_0} $ we see that these two measures satisfy condition a) of Lemma \ref{lemma:metric-ratio-singleton} for $c=1$ and any probability measure $\eta$ on $X$. Let
$\Lambda = \{ \frac{i}{n}: i=0, \ldots, n \}$ for some large $n$. Then the unique set of measures in $M_{\Lambda} (\delta_{x_0},\delta_{x'_0}) $ is precisely the set $\{ \mu_i= (1-\frac{i}{n})\delta_{x_0} + \frac{i}{n}\delta_{x'_0}: i = 0, \ldots, n \} $ since $\Phi$
is an isometry with $\{ \Phi(\delta_{x_0}), \Phi(\delta_{x'_0})\} = \{ \delta_{x_0}, \delta_{x'_0} \} $ it follows 
that $$\{ \Phi(\mu_i): i= 0, \ldots, n \} = \{ \mu_i: i= 0, \ldots, n\}.$$
Letting $n \to \infty$ we obtain that 
$$ \Phi(  \{ (1-t)\delta_{x_0} + t \delta_{x'_0}: t\in [0,1] \}) = \{ (1-t)\delta_{x_0} + t \delta_{x'_0}: t\in [0,1] \}.$$

Since by assumption $X$ contains at least three points, it follows that there must be at least one more point $x_1 \neq x'_0$ such that $(x_0, x_1) \in \mathfrak{F}$. In the same way, there must be at least one more point $x'_1 \neq x_0$ such that $(x'_0, x'_1) \in \mathfrak{F}$. Note, that it may happen that $x_1 = x'_1$. 
Let us consider the (nonempty) sets:

$$ \mathcal{N}:= \mathcal{M} \setminus \{ (1-t)\delta_{x_0} + t \delta_{x'_0}: t\in [0,1] \} \ \text{and} \ \mathcal{N'}:= \mathcal{M'} \setminus \{ (1-t)\delta_{x_0} + t \delta_{x'_0}: t\in [0,1] \}.$$

Then we have that $\Phi(\mathcal{N}) = \mathcal{N'}$. Since $\Phi$ is an isometry with $\Phi(\delta_{x_0}) = \delta_{x'_0} $ we have:
$$ \max \{d_{\Wp}(\delta_{x_0}, \nu) \, | \, \nu \in \mathcal{N} \} =\{d_{\Wp}(\delta_{x'_0}, \nu) \, | \, \nu \in \mathcal{N'} \}. $$
Furthermore, we have 
$$ \max \{d_{\Wp}(\delta_{x_0}, \nu) \, | \, \nu \in \mathcal{N} \} = \max \{d(x_0,x) \,|\, (x_0,x) \in \cF\},x \neq x'_0 \}, $$ 
and 
$$ \max \{d_{\Wp}(\delta_{x'_0}, \nu') \, | \, \nu' \in \mathcal{N'} \} = \max \{d(x'_0,x') \,|\, (x'_0,x') \in \cF\},x' \neq x_0 \}. $$ 

Let us denote by $\hat{x} \in X$ the point where the first and by $\hat{x}' \in X$ the point where the second maximum is achieved. By the above consideration we have:
$$ d(x_0, \hat{x}) = d(x'_0, \hat{x}'),$$
which is in contradiction to the anisometric property of $X$, finishing the proof. 
\end{proof}

\subsection{Proof of Lemma \ref{lemma:metric-ratio-singleton}.}

In this subsection we give the proof of Lemma \ref{lemma:metric-ratio-singleton}.  
\begin{proof}
    Let us start with (a) $\Rightarrow$ (b). Assume that $(v,v') \in \mathfrak{F}$ defined by \eqref{eq:frak-F-def}, $c \in [0,1], \, \eta \in \cPX,$ and $\mu=(1-c) \eta + c \delta_v$, and $\nu=(1-c)\eta + c \delta_{v'}$.

    In order to conclude (b) we start the proof by some preparation. Let us define the function $\varphi_v:\, X \to [0, \infty)$ as follows:
    \begin{align} \label{eq:phi-v-def}
        \varphi_v(y):= \min\left\{\sum_{m=0}^{M-1} d^p(x_m, x_{m+1}) \, \middle| \ M \geq 1 \,
         (v=x_0, x_1, \dots, x_{M-1},x_M=y) \right. \\ \left. \text{is a chain of points from}\  v \ \text{to} \  y \right\}.
    \end{align}
    A crucial property of $\varphi_v$ is that 
    \begin{align} \label{eq:phi-v-is-competitive}
        \varphi_v(y)-\varphi_v(x) \leq d^p(x,y) \text{ for all } x,y \in X.
    \end{align}
    Indeed, \eqref{eq:phi-v-is-competitive} is equivalent to $\varphi_v(y) \leq \varphi_v(x)+d^p(x,y),$ which is easily shown by the following: if $(v=x_0, x_1, \dots, x_{K-1},x_K=x)$ is a chain connecting $v$ and $x$ such that $\sum_{k=0}^{K-1} d^p(x_k,x_{k+1})=\varphi_v(x),$ then the chain $(v=x_0, x_1, \dots x_{K-1}, x_K=x, x_{K+1}=y)$ connects $v$ and $y,$ and $$\sum_{k=0}^K d^p(x_k, x_{k+1})=\sum_{k=0}^{K-1} d^p(x_k,x_{k+1}) + d^p(x_K,x_{K+1})=\varphi_v(x)+d^p(x,y),$$ and hence $\varphi_v(y) \leq \varphi_v(x)+d^p(x,y)$ as $\varphi_v(y)$ is the minimum taken over all chains connecting $v$ and $y.$
    \par
    
    A consequence of \eqref{eq:phi-v-is-competitive} is that 
    \begin{align} \label{eq:weak-Kant-duality}
        \int_X \varphi_v ~\dd (\sigma -\rho) \leq d_{\cW_p}^p(\rho, \sigma)
    \end{align}
    holds for all $\rho, \sigma \in \cPX.$ Indeed, \eqref{eq:weak-Kant-duality} is an instance of the weak Kantorovich duality, and follows from the direct computation
    \begin{align}
        \int_X \varphi_v ~\dd (\sigma -\rho)&=\int_X \varphi_v(y) ~\dd \sigma(y)-\int_X \varphi_v(x) ~\dd \rho(x) \nonumber \\
        & =\iint_{X \times X} \ler{\varphi_v(y) -\varphi_v(x)} ~\dd \pi_{\rho, \sigma}(x,y) \nonumber\\
        &\leq \iint_{X \times X} d^p(x,y) ~\dd \pi_{\rho, \sigma}(x,y)
        =d_{\cW_p}^p(\rho, \sigma),
    \end{align}
    where $\pi_{\rho, \sigma}$ is an optimal coupling of $\rho$ and $\sigma.$
    
    Moreover, in fact we have that $\int_X \varphi_v ~\dd (\nu-\mu)=d_{\cW_p}^p(\mu,\nu)$. This is, because
    \begin{align}
        \int_X \varphi_v ~\dd (\nu-\mu)
        =\int_X \varphi_v ~\dd (c\ler{\delta_{v'}-\delta_v})
        =c \ler{\varphi_v(v')- \varphi_v(v)}
        =c d^p(v,v'),
        \nonumber
    \end{align}
    where we used that $\varphi_v(v)=0$ and $(v,v') \in \mathfrak{F}$ and hence $\varphi_v(v')=d^p(v,v').$ 
    
    Furthermore, note, that  the Kantorovich cost of the coupling $$\pi_{c}:=(1-c)(\mathrm{Id}_X \times \mathrm{Id}_X)_{\#}\eta +c \delta_{(v,v')}$$ is $K(\pi_c)=c d^p(v,v'),$ which implies that $\pi_c$ is optimal and $d_{\cW_p}^p(\mu,\nu)=c d^p(v,v').$
    \par
    It follows from \eqref{eq:weak-Kant-duality} that for any $\xi_1, \dots, \xi_N \in \cPX$ the inequality 
    \begin{align} \label{eq:lower-bound-on-sum-of-costs}
    d_{\cW_p}^p(\mu,\xi_1)+d_{\cW_p}^p(\xi_1,\xi_2)+ \dots, d_{\cW_p}^p(\xi_{N-1},\xi_N)+d_{\cW_p}^p(\xi_N,\nu) \geq d_{\cW_p}^p(\mu,\nu)
    \end{align}
    holds. To justify \eqref{eq:lower-bound-on-sum-of-costs}, note that $d_{\cW_p}^p(\xi_j,\xi_{j+1}) \geq \int_X \varphi_v ~\dd \ler{\xi_{j+1}-\xi_j}$ by \eqref{eq:weak-Kant-duality} for all $j \in \{0, \dots, N \},$ and hence 
    \begin{align}
        \sum_{j=0}^N d_{\cW_p}^p(\xi_j,\xi_{j+1}) 
        &\geq \sum_{j=0}^N \int_X \varphi_v ~\dd \ler{\xi_{j+1}-\xi_j}
        = \int_X \varphi_v ~\dd \ler{\sum_{j=0}^N\ler{\xi_{j+1}-\xi_j}}
        \nonumber \\
        &= \int_X \varphi_v ~\dd \ler{\xi_{N+1}-\xi_1}
        = \int_X \varphi_v ~\dd \ler{\nu-\mu}
        =d_{\cW_p}^p(\mu, \nu).
    \end{align}
    
    After this preparation we come back to the proof of the implication (a) $\Rightarrow$ (b). Assume that $\lers{\xi_1, \dots, \xi_N}$ belongs to the multiple ratio set $M_{\Lambda}(\mu, \nu)$ for some collection of weights $\Lambda=(\lambda_i)_{i=1}^N,$ then by the definition of $M_\Lambda(\mu, \nu)$ (see \eqref{eq:multiple-ratio-set-def}), the reverse inequality 
    \begin{align}
        \sum_{j=0}^{N} d_{\cW_p}^p(\xi_j, \xi_{j+1})
        \leq \ler{\sum_{j=0}^{N} \ler{\lambda_{j+1}-\lambda_j}}d_{\cW_p}^p(\mu, \nu)=d_{\cW_p}^p(\mu, \nu)
    \end{align}
    holds. But this means that $d_{\cW_p}^p(\xi_j,\xi_{j+1})=\int_X \varphi_v ~\dd \ler{\xi_{j+1}-\xi_j}$ for all $j \in \{0, \dots, N\}.$ 
    
    Furthermore, the equality $d_{\cW_p}^p(\xi_j,\xi_{j+1})=\int_X \varphi_v ~\dd \ler{\xi_{j+1}-\xi_j}$ implies that for all optimal coupling $\pi_{j,j+1}$ between $\xi_j$ and $\xi_{j+1},$ the equality 
    \begin{align} \label{eq:tightness}
        \varphi_v(y)-\varphi_v(x)=d^p(x,y)
    \end{align}
    holds whenever $\pi_{j,j+1}\ler{\lers{(x,y)}}>0.$ Indeed, the inequality
    \begin{align}
        \int_X \varphi_v ~\dd \ler{\xi_{j+1}-\xi_j}
        &=\iint_{X \times X} \ler{\varphi_v(y)-\varphi_v(y)} ~\dd \pi_{j,j+1}(x,y) \nonumber\\ 
        &\leq
        \iint_{X \times X} d^p(x,y) ~\dd \pi_{j,j+1}(x,y)=d_{\cW_p}^p(\xi_j,\xi_{j+1})
    \end{align}
    is saturated only if \eqref{eq:tightness} holds on the full support of $\pi_{j,j+1}.$
    
    In what follows we shall study the structure of the measures $\xi_j$ by investigating the supports of $\pi_{j, j+1}$ that describes the way in which mass is transported from $\xi_j$ to $\xi_{j+1}$ step by step. Our special interest is to look at the support points in $v$ and $v'$. 
    
    Let us  note first, that \eqref{eq:tightness} cannot hold if $y=v$ and $x \neq v$ because the left-hand side of \eqref{eq:tightness} is negative in this case while the right-hand side of \eqref{eq:tightness} is positive. Moreover \eqref{eq:tightness} cannot hold even if $y=v'$ and $x \in X\setminus \{v,v'\}.$ Indeed, for an arbitrary $x \in X\setminus \{v,v'\},$ let $(v=x_0, x_1, \dots, x_M=x)$ be a chain from $v$ to $x$ such that $\sum_{m=0}^{M-1}d^p(x_m, x_{m+1})=\varphi_v(x).$ Then by the assumption $(v, v') \in \mathfrak{F}$ one gets that for the path $(v=x_0,x_1, \dots, x_M=x, x_{M+1}=v')$ connecting $v$ and $v'$ the strict inequality  
    \begin{align}
        \varphi_v(v')=d^p(v,v')<\sum_{m=0}^M d^p(x_m, x_{m+1})=\sum_{m=0}^{M-1} d^p(x_m, x_{m+1})+d^p(x,v')=\varphi_v(x)+d^p(x,v')
        \nonumber
    \end{align} 
    holds. 
    
    Let us introduce the partitions $P_1:=\lers{v,v'}$ and $P_2:=X \setminus \{v,v'\}$ of $X.$ We just proved that \eqref{eq:tightness} cannot hold if $y \in P_1$ and $x \in P_2,$ which in turn implies that for all $j \in \{0, \dots, N\}$ and for all $\pi_{j,j+1}$ optimal coupling between $\xi_j$ and $\xi_{j+1},$ one has
    $\pi_{j,j+1}(P_2 \times P_1)=0.$ That is, no mass is transported from the partition $P_2$ to $P_1,$ which implies that $\xi_{j+1}(P_2) \geq \xi_j(P_2)$ for all $j.$ 
    
    Therefore we have that:
    \begin{align}
        \mu(P_2) =\xi_0(P_2)\leq \xi_1(P_2) \leq \dots \leq \xi_N(P_2) \leq \xi_{N+1}(P_2)=\nu(P_2).
        \nonumber
    \end{align}
    On the other hand, from the form of $\mu$ and $\nu$  we know that 
    \begin{align}
    \nu(P_2)-\mu(P_2)=\ler{c\ler{\delta_{v'}-\delta_v}}(P_2)=0,    \nonumber
    \end{align}
    and hence $\xi_{j+1}(P_2)=\xi_j(P_2)$ for all $j \in \{0,\dots, N\}.$ But then we can conclude that no mass is transported from $P_1$ to $P_2$ either at any step, that is, $\pi_{j, j+1}(P_1 \times P_2)=0$ for all $j$ and optimal $\pi_{j, j+1},$ because $\pi_{j, j+1}(P_1 \times P_2)>0$ would contradict that $\xi_{j+1}(P_2)=\xi_j(P_2)$ and $\pi_{j, j+1}(P_2 \times P_1)=0.$ So $\pi_{j, j+1}(P_1 \times P_2)=\pi_{j, j+1}(P_2 \times P_1)=0$ for all $j$ and optimal $\pi_{j, j+1},$ and we also noted that $\pi_{j,j+1}(\{(v',v)\})=0,$ which implies that 
    \begin{align} \label{eq:diff-of-weights-by-pi-j-j+1}
        \xi_{j+1}(\{v'\})-\xi_j(\{v'\})=\pi_{j,j+1}(\{(v,v')\}).
    \end{align}
Furthermore, by \eqref{eq:diff-of-weights-by-pi-j-j+1} we get that
\begin{align} \label{eq:decomp-of-costs}
    c d^p(v,v')=d_{\cW_p}^p(\mu, \nu)
    =\sum_{j=0}^N \iint_{X^2} d^p(x,y) ~\dd \pi_{j,j+1}(x,y) \nonumber \\
    =\sum_{j=0}^N \ler{\pi_{j,j+1}(\{(v,v')\}) d^p(v,v') + \iint_{X^2\setminus \{(v,v')\}} d^p(x,y) ~\dd \pi_{j,j+1}(x,y)} \nonumber \\
    =\sum_{j=0}^N \ler{\xi_{j+1}(\{v'\})-\xi_j(\{v'\})} d^p(v,v')
    +\sum_{j=0}^N \iint_{X^2\setminus \{(v,v')\}} d^p(x,y) ~\dd \pi_{j,j+1}(x,y).
\end{align}
The first term of the last line of \eqref{eq:decomp-of-costs} is a telescopic sum that can be computed as
\begin{align}
    \sum_{j=0}^N \ler{\xi_{j+1}(\{v'\})-\xi_j(\{v'\})} d^p(v,v')
    &=\ler{\xi_{N+1}(\{v'\})-\xi_0(\{v'\})} d^p(v,v')
    \nonumber \\ 
    &=(\nu-\mu)(\{v'\}) d^p(v,v')=c d^p(v,v'),
    \nonumber
\end{align}
which implies that the second term of that line vanishes:
\begin{align} \label{eq:second-term-vanishes}
\sum_{j=0}^N \iint_{X^2\setminus \{(v,v')\}} d^p(x,y) ~\dd \pi_{j,j+1}(x,y)=0.
\end{align}
As all the distances between distinct points of $X$ are positive, \eqref{eq:second-term-vanishes} implies that 
\begin{align} \label{eq:pi-j-j+1-is-almost-trivial}
    \pi_{j,j+1}=(1-c_j)(\mathrm{Id}_X \times \mathrm{Id}_X)_{\#}\eta_j + c_j \delta_{(v,v')}
\end{align}
for some $\eta_j \in \cPX$ and $c_j=\xi_{j+1}(\{v'\})-\xi_j(\{v'\}).$ But then
\begin{align}
    \xi_{j+1}=\xi_j+c_j (\delta_{v'}-\delta_v),
\end{align}
which implies that 
\begin{align}
    \xi_J
    &=\xi_0+\sum_{j=0}^{J-1} c_j (\delta_{v'}-\delta_v)
    =\mu +\ler{\sum_{j=0}^{J-1} \xi_{j+1}(\{v'\})-\xi_j(\{v'\})} (\delta_{v'}-\delta_v) \nonumber \\
    &=\mu +\ler{\xi_J(\{v'\})-\mu(\{v'\})} (\delta_{v'}-\delta_v).
\end{align}
Consequently, there exist weights $\lambda_1, \dots, \lambda_N$ such that $\xi_J=(1-\lambda_J) \mu+\lambda_J \nu$ for all $J \in \{1, \dots, N\},$ as desired, and this concludes the proof of the direction (a) $\Rightarrow$ (b). 
Note that here we used the inequality $\ler{\xi_J(\{v'\})-\mu(\{v'\})} \geq 0$ which follows from \eqref{eq:diff-of-weights-by-pi-j-j+1}. 

\par

We turn now to the direction $(b)$ implies $(a)$. To see this, assume first by contradiction that $(b)$ holds with $\mu,\nu$ given by $\mu=(1-c)\eta + c \delta_v$, $\nu=(1-c)\eta  + c \delta_{v'}$, but such that $(v,v') \notin \mathfrak{F}$. In particular, there exists a chain of (at least three) points $\{x_i\}_{i=0}^{N+1}$ with $x_0=v,x_{N+1}=v'$ such that 
$$
d^p(v, v') \geq \sum_{i=1}^{N}d^p(x_i, x_{i+1}).
$$

Consider the optimal transport plan $\pi_0$ between $\mu$ and $\nu$. If there exists $w \in X$, $w \neq v, v'$ such that $P:=\pi_0(v, w)>0$, we consider the transport plans
\begin{align*}
\pi_1 =& \pi_0 + P(\delta_{(v, v)} - \delta_{(v, w)}), \\
\pi_2 =& (\mathrm{Id}_X \times \mathrm{Id}_X )_{\#}\xi + P(\delta_{(v, w)} - \delta_{(v, v)}),
\end{align*}
with $\xi=p_{2\#}\pi_1$. Then $\pi_1$ is a transport plan between $\mu$ and $\xi$, $\pi_2$ is a transport plan between $\xi$ and $\nu$ such that 
$$
K(\pi_1) + K(\pi_2)=K(\pi_0).$$

\noindent We claim, that if $\sum_{i=1}^{N+1}K(\pi_i) \leq K(\pi_0)$ then there exists $\Lambda\in [0,1]^N$ such that $\{\xi_i\}_{i=1}^N \in M_{\Lambda}(\mu, \nu)$. To see this, set  $\lambda_i=(\sum_{j=1}^{i}K(\pi_{j}))/K(\pi_0)$ for $i \leq N$. Since costs are positive it is clear that $\Lambda = \{\lambda_i\}_{i=1}^N \in [0,1]^N$ satisfies the monotonicity condition. Then we have that for $i \leq N$,  
$$
d_{\mathcal{W}_p}^p(\xi_{i-1}, \xi_{i}) \leq K(\pi_{i})=\sum_{j=1}^{i}K(\pi_j)-\sum_{j=1}^{i-1}K(\pi_j)=(\lambda_{i}-\lambda_{i-1})d^p_{\mathcal{W}_p}(\mu, \nu),
$$ as claimed. Therefore, $K(\pi_1) + K(\pi_2)=K(\pi_0)$ implies that $\xi\in M_{\lambda_0}(\mu, \nu)$ for some $\lambda_0$; since $\xi(\{w\})>\mu(\{w\})$ and $w \neq v,v'$, we thus have that $\xi\neq (1-\lambda_0)\mu + \lambda_0 \nu$, giving the contradiction that $M_{\lambda_0}(\mu, \nu)$ is not a singleton. 

If $\pi_0(v, w)=0$ for any $w \neq v,v'$, then in particular $\pi_0(v, v')\geq c$, and $K(\pi_0)\geq cd^p(v, v')$. Consider now the measures 
$$
\xi_i=(1-c)\eta + c \delta_{x_i}
$$ for $1 \leq i \leq N$, with $x_i$ given by the chain above. Then 
$$
\pi_{i}=(1-c)(\mathrm{Id}_X \times \mathrm{Id}_X)_{\#}\eta + c \delta_{(x_{i-1}, x_i)} 
$$ is a transport plan between $\xi_{i-1}$ and $\xi_i$ with Kantorovich cost $K(\pi_i)=cd^p(x_{i-1}, x_1)$.

Since by the assumption that $(v, v') \notin \mathfrak{F}$ we have that $\sum_{i=1}^{N+1}K(\pi_i) \leq K(\pi_0)$. Again, by the claim above,  there exists $\Lambda\in [0,1]^N$ such that $\{\xi_i\}_{i=1}^N \in M_{\Lambda}(\mu, \nu)$.

Now, by their definition it is clear that the measures $\xi_i$ are not obtained by a simple interpolation between $\mu$ and $\nu$, which gives the needed contradiction.  

We now consider the second case, that is we assume by contradiction that there exists no measure $\eta$ and no weight $c$ such that $\mu=(1-c)\eta + c \delta_v$, $\nu=(1-c)\eta + c \delta_{v'}$. We claim that if this holds, then there must exist  two points $x_1 \neq x_2$ such that either $\mu(x_1) > \nu(x_1)$, $\mu(x_2) >\nu(x_2)$ or $\nu(x_1) > \mu(x_1)$, $\nu(x_2) >\mu(x_2)$. 
We shall prove the claim by contradiction. If no such pair of two different points exist then the first possibility is that the two measures are equal: $\mu=\nu$. But this will lead to a contradiction by taking $c=0$ and $\eta := \mu = \nu$. 
The other possibility, when the measures are not equal, this will imply that there exist  exactly one pair of points $(x_1, x_2)$ where $\mu$ and $\nu$ differ such a way that $\mu(x_1)> \nu(x_1)$ and $\nu(x_2) > \mu(x_2)$. In this case we define $\eta:=\delta_{x_1}$ if $\mu=\delta_{x_1}$ and
$$ \eta:= \frac{\mu-\mu(x_1)\delta_{x_1}}{1-\mu(x_1)}$$ otherwise. Then clearly we have $\mu=(1-\mu(x_1))\eta + \mu(x_1)\delta_{x_1}$. Moreover, since the two measures differ only in the points $x_1$ and $x_2$ we must have that $\nu=(1-\mu(x_1))\eta + \mu(x_1)\delta_{x_2}$. This representation of the two measures is in contradiction to the assumption finishing the proof of our claim.

We assume without loss of generality that $\mu(x_1) < \nu(x_1)$, $\mu(x_2) <\nu(x_2)$. Then, if $\pi_0$ is the optimal coupling between $\mu$ and $\nu$, and $x_1'\in X$ is such that $\pi_0(x_1', x_1)>0$ (which must exist, as mass is transferred from somewhere else to $x_1$), we consider the measure $\xi \in \cPX$ given by 
$$
\xi=\nu +\pi_0(x_1', x_1)(\delta_{x_1'}-\delta_{x_1}). 
$$ 
Observe that the measure on the product space given by 
$$
\pi_1=\pi_0 + \pi_0(x_1',x_1)(\delta_{(x_1', x_1')}-\delta_{(x_1', x_1)})
$$ is a coupling between $\mu$ and $\xi$. Similarly, the formula 
$$
\pi_2=(1-\pi_0(x_1',x_1)) (\mathrm{Id}_X \times \mathrm{Id}_X)_{\#}\xi + \pi_0(x_1',x_1)\delta_{(x_1', x_1)}
$$ defines a coupling between $\xi$ and $\nu$. Due to its construction, $\xi$ cannot be written as $\xi=(1-\lambda)\mu + \lambda\nu$. Indeed, while $\mu(x_1)\leq \xi(x_1)<\nu(x_1)$, we have $\mu(x_2)<\xi(x_2)=\nu(x_2)$. However, we also have that $K(\pi_1) + K(\pi_2)=K(\pi_0)$, and thus $\xi \in M_{\lambda}(\mu, \nu)$ for some $\lambda$ as we have seen before. In particular, the set $M_{\lambda}(\mu, \nu)$ is not a singleton. This contradiction finishes the proof of the fact that $(b) $ implies $(a)$.
\end{proof}

\section{Spaces with rigid $p$-Wasserstein spaces form a $G_{\delta}$ subset of all compact metric spaces} 

Let us recall that in a topological space a $G_{\delta}$ set is a countable intersection of open sets. The main result of this section is the following:
\begin{proposition} \label{prop:rigid-G-delta}
 For any parameter $p \geq 1,$ the subset $\mathcal{R}_p\subset \mathfrak{M}$ of compact metric spaces $(X, d)$ with the property that $\mathrm{Isom(}\Wp(X, d))= \{\mathrm{Id}_{\Wp(X,d)} \}$ is a $G_{\delta}$ set in the spaces of compact metric spaces $\mathfrak{M}$ with respect to the Gromov-Hausdorff topology.       
\end{proposition} 

The proof of Proposition \ref{prop:rigid-G-delta} is based on Lemma \ref{lemma:R-G-delta} below. Before going into the details, let us recall that if $(X,d)$ is compact, then $\Wp(X, d)$ is also compact (see e.g. \cite[Section 6]{Villani}). Therefore, we can consider for $p\geq 1$ the Wasserstein map
$$ T_p: \mathfrak{M} \to \mathfrak{M} \ \text{defined as:} \  T_p((X, d)) = (\Wp(X, d), d_{\mathcal{W}_p}).$$
Note, that if $(X, d_X)$ and $(Y, d_Y)$ are isometric then the same is true for the associated Wasserstein spaces. Indeed, if $\varphi: (X, d_X) \to (Y, d_Y)$ is an isometry. Then the push-froward map $\Phi = \varphi_{\sharp}$ defines an isometry between $(\Wp(X, d_X), d_{\mathcal{W}_p})$ and $(\Wp(Y, d_Y), d_{\mathcal{W}_p})$. Therefore the map $T_p$ can be understood as being defined between equivalence classes of metric spaces induced by isometries. It is known that this Wasserstein map is $1-$ Lipschitz with respect to the Gromov-Hausdorff metric:
\begin{equation} \label{eq:1-lipshitz}
d_{GH}(T_p(X), T_p(Y)) \leq d_{GH}(X, Y). 
\end{equation}
For the proof see \cite[Theorem 3.13]{Memoli}.  In the following let us denote by 
$$ \mathcal{R}:= \{ X \in \mathfrak{M}: \ \mathrm{Isom}(X) = \{ \mathrm{Id}_X \} \}.$$ The following lemma shows that $\mathcal{R}$ is a $G_{\delta}$ set in $\mathfrak{M}$.

\begin{lemma} \label{lemma:R-G-delta}
For every $n\geq 1$ the set 
$$ \mathcal{R}_n := \{ X \in \mathfrak{M}: \max_{x \in X} d(x, F(x)) < \frac{1}{n} \ \text{for all} \ F\in \mathrm{Isom}(X) \} $$
is an open set in $\mathfrak{M}$. Consequently, the set 
$\mathcal{R}= \bigcap_{i=1}^{\infty} R_n$ is a $G_{\delta}$ set in $\mathfrak{M}$.
\end{lemma}

    \begin{proof}
  More generally we shall show that for all $\varepsilon >0$  the set  
  $$ \mathcal{R}_{\varepsilon} := \{ X \in \mathfrak{M}: \max_{x \in X} d(x, F(x)) < \varepsilon \ \text{for all} \ F\in \mathrm{Isom}(X) \} $$
  is open in $\mathfrak{M}$.  We shall do that by showing that its complement in $\mathfrak{M}$:
$$ \mathfrak{M}\setminus \mathcal{R}_{\varepsilon}= \mathcal{Q}_{\varepsilon} =  \{ X \in \mathfrak{M}: \max_{x \in X} d(x, F(x)) \geq \varepsilon \ \text{for some} \ F\in \mathrm{Isom}(X) \} $$ 
is closed in $\mathfrak{M}$. 

To do so let us choose a sequence $(X_k)_{k\in\mathbb{N}} \subset \mathcal{Q}_{\varepsilon}$ such that 
$d_{GH}(X_k, X) \to 0$ as $k\to \infty$ for some $X\in \mathfrak{M}$. We have to show that 
$X\in \mathcal{Q}_{\varepsilon}$. 

By passing to isometric copies we may assume that $X_k \subset Z$ for all $k$ and also 
$X\subset Z$ for some $Z\in \mathfrak{M}$. Since $X_k \in \mathcal{Q}_{\varepsilon}$ or each $k$,  there exists $F_k\in \mathrm{Isom}(X_k)$ and $x_k \in X_k$ with the property that $d(x_k, F_k(x)) \geq \varepsilon$. 

Our aim to show that there exists $F_0\in \mathrm{Isom}(X)$ and $x_0\in X$ such that $d(x_0, F_0(x_0)) \geq \varepsilon$. The map $F_0$ will be obtained as the limit of some subsequence  of $(F_k)_{k\in\mathbb{N}}$. Similarly the point $x_0$ will also be obtained as the limit by a subsequence of $(x_k)_{k\in\mathbb{N}}$. In order to explain the construction of $F_0$ let us choose a point $x \in X$. Since $d_H(X_k, X) \to 0$ we find a sequence $(x_k)_{k\in\mathbb{N}}$ such that $x_k \in X_k$ and $x_k \to x$. 
Now let us consider the sequence of image points $(F_k(x_k))_{k\in\mathbb{N}} \subset Z$. By compactness we find a  subsequence $(x_{k_j})_{j\in\mathbb{N}}$ such that  $(F_{k_j}(x_{k_j}))_{j\in\mathbb{N}}$ converges in $Z$. It is easy to see that that the limit point is in $X$. Let us call this limit point $F_0(x)$. Now try to repeat this argument for all points $x\in X$. The emerging issue is that the converging subsequence might depend on the point. The solution to this problem is as follows: notice first that since $Z$ is a compact metric space, it is separable and therefore $X$ is separable as well. 
Let $Y\subset X$ be a countable dense subset of $X$. By a diagonal argument we can obtain a sequence $(y_k)_{k\in\mathbb{N}}$ such that for all $y \in Y$ there exists a subsequence $(y_{k_j})_{j\in\mathbb{N}}$ such that
$y_{k_j} \in X_{k_j}$, $y_{k_j} \to y$ and $F_{k_j}(y_{k_j})$ converges to some point in $X$. By this procedure we can define a map $F_0: Y \to X$. By passing to this subsequence we can assume without loss of generality that for $y \in Y$ we can find a sequence $y_k \in X_k$ such that $y_k \to y$ and  $F_k(y_k) \to F_0(y) $. 

Let us notice first that this construction is independent on the choice of $(y_k)_{k\in\mathbb{N}}$. Indeed, if we choose another sequence  $(y'_k)_{k\in\mathbb{N}} $ with $y'_k\in X_k$ such that $y'_k\to y$
then since $F_k$ is an isometry of $X_k$ we have $d(F_k(y_k), F_k(y'_k)) = d(y_k, y'_k) \to 0$ as $k \to \infty$ and thus $F_k(y'_k) \to F_0(y)$ as well. 

The second observation is that the map $F_0$ is an isometric embedding of $Y$ into $X$. 
Indeed, let  $y, y' \in Y$ and $y_k \to y$ $y'_k \to y'$ be the corresponding approximating sequences. Since $F_k$ is an isometry we have 
$$ d(F_k(y_k), F_k(y'_k)= d(y_k, y'_k),$$
and passing to the limit as $k \to \infty$ we obtain
$$ d(F_0(y), F_0(y')) = d(y, y') \ \text{for any} \ y, y' \in Y.$$
Using this property we can extend $F_0 $ to a map on the whole of $X$ obtaining an isometric embedding $F_0: X \to X$. It turns out that this mapping is also surjective since the approximating isometries $F_k: X_k \to X_k$ were surjective and so we have that 
$F_0 \in \mathrm{Isom}(X)$. 

To finish the proof we shall consider the sequence of points $x_k \in X_k $ having the property that 
$$d(F_k(x_k), F_k(x_k) ) \geq \varepsilon.$$
By compactness of $Z$ we find a converging subsequence (denoted in the same way) $x_k \to x_0$, where clearly $x_0 \in X$. Passing to the limit in the above inequality as $k\to \infty $ we obtain 
$$d(F_0(x_0), F_0(x_0) \geq \varepsilon,$$
which shows that $X \in \mathcal{Q}_{\varepsilon}$, concluding the proof. 
    \end{proof}

Combining the above results the proof of Proposition \ref{prop:rigid-G-delta} is immediate:

%\begin{proof} {\it (Proposition \ref{prop:rigid-G-delta})}. 
\begin{proof}[Proof of Proposition \ref{prop:rigid-G-delta}]
Consider the map 
$$T_p:\mathfrak{M} \to \mathfrak{M} ,  \ X \to \Wp(X,d).$$
 By  \eqref{eq:1-lipshitz} this map is Lipschitz continuous. By Lemma \ref{lemma:R-G-delta} we obtain that the set 
 $$ \mathcal{R}:= \{ X \in \mathfrak{M}: \ \mathrm{Isom}(X) = \{ \mathrm{Id}_X \} \}$$ 
 is a $G_{\delta}$  subset of $\mathfrak{M}$. Since preimages of $G_{\delta}$  sets under continuous maps are again $G_{\delta}$ sets we obtain that $T_p^{-1}(\mathcal{R})$ is also a $G_{\delta}$ set. 
 But on the other hand
 $$T_p^{-1}(\mathcal{R}) = \{ X \in \mathfrak{M}: \mathrm{Isom}(\Wp(X,d) )= \{ \mathrm{Id}_X \} \}$$ proving the statement of the Proposition. 
\end{proof}

Let us give a short summary of the proof of Theorem \ref{thm:rigidity-is-generic} that is already indicated in the introduction. The proof is based on Corollary \ref{rigid-dense} and Proposition \ref{prop:rigid-G-delta}. Note first, that the set 
\begin{equation} \label{eq:R_p}
  \mathcal{R}_p = \{ X \in \mathfrak{M}: \ \mathrm{Isom}(\Wp(X,d)) = \{ \mathrm{Id}_X \} \}  
\end{equation} 
is a subset of metric spaces $X \in \mathfrak{M}$ for which the $\Wp(X,d)$ is rigid. On the other hand, Corollary \ref{rigid-dense} and Proposition \ref{prop:rigid-G-delta} show that $ \mathcal{R}_p$ is a dense $G_{\delta}$ set in $\mathfrak{M}$ proving Theorem \ref{thm:rigidity-is-generic}.

\section{The Dirac preservation property and the Wasserstein recovery property is generic}

An interesting open question concerning Wasserstein spaces is whether the isometry of $\Wp(X,d_X)$ and $\Wp(Y,d_Y)$ implies that the underlying metric spaces $(X,d_X)$ and $(Y,d_Y)$ are also isometric. We cannot answer this question in full generality, but we can answer it for generic objects. Since we will occasionally consider different metrics on the same underlying set in this section, we will use slightly more detailed notation. In order to formulate our result we start with two definitions:

\begin{definition}[Recovery property] \label{def:RECp}
Let $p\geq 1.$ We say that the Polish metric space $(X, d_X)$ satisfies the \emph{recovery property of order $p$ (RECp)} if for any $(Y, d_Y)$ such that $\Wp(X,d_X)$ is isometric to $\Wp(Y,d_Y)$
we have that $X$ is isometric to $Y$.
\end{definition} 

\begin{definition}[Dirac preservation property] \label{def:DPPp}
Let $p\geq 1.$ We say that the Polish metric space $(X, d_X)$ satisfies the \emph{Dirac preservation property of order $p$ (DPPp)} if for any $(Y, d_Y)$ Polish metric space and any isometry $\Phi: \Wp(X,d_X) \to \Wp(Y,d_Y),$ the measure $\Phi(\delta_x) \in \Wp(Y,d_Y)$ is a Dirac mass for all $x \in X,$ and for all $y \in Y$ there exists some $x \in X$ such that $\Phi(\delta_x)=\delta_y.$ 
\end{definition}

Let us note first that the Dirac preservation property implies the recovery property. Indeed, let $(X,d_X)$ satisfy the DPPp for some $p \geq 1,$ and assume that $\Wp(Y,d_Y)$ is isometric to $\Wp(X,d_X)$ for some Polish metric space $(Y, d_Y).$ Let $\Phi: \Wp(X,d_X) \to \Wp(Y,d_Y)$ be an isometry. By the DPPp, it maps Dirac masses to Dirac masses, that is, for every $x\in X$ there exists a unique $y \in Y$ such that $\Phi(\delta_x) = \delta_y.$ Accordingly, let us define the map $F:X \to Y$ by the equation $\Phi(\delta_x) = \delta_{F(x)}.$
Then, $F$ is an isometric embedding of $X$ into $Y$. 
Indeed, for any $x,x' \in X$ one gets that 
\begin{align}
d_Y(F(x), F(x')) = d_{\Wp(Y,d_Y)}(\delta_{F(x)} , \delta_{F(x')}) &= d_{\Wp(Y,d_Y)}(\Phi(\delta_{x}), \Phi(\delta_{x'})) \nonumber \\ 
&=d_{\Wp(X,d_X)}(\delta_x, \delta_{x'} ) = d_X(x, x').
\end{align}
Moreover, as DPPp ensures that for every Dirac mass $\delta_y \in \Wp(Y,d_Y)$ there is an $x\in X$
such that $\Phi(\delta_x) = \delta_y,$ we obtain that the map $F:X \to Y$ is surjective and thus an isometry from $X$ to $Y.$ In particular, $X$ and $Y$ are isometric, and hence $(X,d_X)$ satisfies the recovery property RECp. 
\par

It is an interesting question whether the recovery property and the Dirac preservation property are actually equivalent or not. 
The goal of the following argument is to demonstrate that the DPPp property of metric spaces is strictly stronger than RECp property. We will show that for every $p\geq 1,$ the compact metric space $([0,1],|.|^{1/p})$ satisfies the recovery property of order $p$ (RECp) but does not satisfy the Dirac preservation property of order $p$ (DPPp). The shorthand $|.|^{1/p}$ stands for the metric $d(x,y)=|x-y|^{1/p},$ and it is known that $\cW_p([0,1],|.|^{1/p})$ admits mass-splitting isometries \cite{GTV1, GTV2} for all $p\geq 1$, and hence $([0,1],|.|^{1/p})$ does not satisfy the DPPp. 
\par
The next task is to show that $([0,1],|.|^{1/p})$ satisfies RECp, that is, $(Y,d_Y)$ is isometric to $([0,1],|.|^{1/p})$ whenever $\cW_p(Y,d_Y)$ is isometric to $\cW_p([0,1],|.|^{1/p}).$
\par 
The first observation is that 
\begin{align} \label{eq:Wass-correspondence}
    d_{\cW_p ([0,1],|.|^{1/p})}^{p}=d_{\cW_1 ([0,1],|.|)}.
\end{align}
Indeed, for any probability measure $\mu, \nu$ on $[0,1]$ 
\begin{align}
    d_{\cW_p ([0,1],|.|^{1/p})}^{p}(\mu, \nu)
    =\inf_{\pi \in \cC(\mu, \nu)} \iint_{[0,1]^2} |x-y|^{\frac{1}{p}p} ~\dd \pi(x,y)
    =d_{\cW_1 ([0,1],|.|)}(\mu, \nu). \nonumber
\end{align}
The right-hand side of \eqref{eq:Wass-correspondence} is a genuine metric, and hence $d_{\cW_p ([0,1],|.|^{1/p})}^{p}$ satisfies the triangle inequality. Therefore, by the assumption $\cW_p(Y,d_Y) \simeq \cW_p([0,1],|.|^{1/p}),$ the quantity $d_{\cW_p(Y,d_Y)}^{p}$ also satisfies the triangle inequality. Restricting this inequality to Dirac measures shows that 
\begin{align}
    r(x,z):=d_Y^p(x,z) \nonumber
\end{align}
is a genuine metric on $Y.$ Furthermore, 
\begin{align} \label{eq:Wass-correspondence-2}
    d_{\cW_p(Y,d_Y)}^p(\mu, \nu)
    &=\inf_{\pi \in \cC(\mu, \nu)} \iint_{Y\times Y} d^p_Y(x,z) ~\dd \pi(x,z) \nonumber \\
    &=\inf_{\pi \in \cC(\mu, \nu)} \iint_{Y\times Y} r(x,z) ~\dd \pi(x,z)
    =d_{\cW_1 (Y,r)}(\mu, \nu),
\end{align}
which implies that $\cW_1(Y,r)$ is isometric to $\cW_1([0,1], |.|).$ Indeed, let $\Phi$ be an isometry between $\cW_p([0,1],|.|^{1/p})$ and $\cW_p(Y,d_Y)$  that exists by assumption. Then for $\mu,\nu \in \Wp(Y,d_Y)$
\begin{align}
    d_{\cW_1([0,1],|.|)}(\mu, \nu)
    &=d_{\cW_p([0,1],|.|^{1/p})}^p(\mu,\nu) \nonumber \\
    &=d_{\cW_p(Y,d_Y)}^p(\Phi(\mu), \Phi(\nu))=d_{\cW_1(Y,r)}(\Phi(\mu), \Phi(\nu))
\end{align}
holds true by \eqref{eq:Wass-correspondence} and \eqref{eq:Wass-correspondence-2}. So 
\begin{align} \label{eq:W-1-Y-r-isom-W-1-01-std}
 \cW_1(Y,r) \simeq \cW_1([0,1], |.|)   
\end{align}
where $\cW_1([0,1], |.|)$ is compact, which implies that so is $\cW_1(Y,r).$ The Dirac measures in $\cW_1(Y,r)$ form a closed isometric copy of $(Y,r),$ and hence $(Y,r)$ is also compact. Furthermore, 
$$
\mathrm{diam}(Y,r)=\mathrm{diam}(\cW_1(Y,r))=\mathrm{diam}(\cW_1([0,1],|.|))=1.
$$ 
The next observation is that there are no two distinct points of $\cW_1([0,1], |.|)$ that are joined by a unique geodesic. One way to see this is to recall that the $1$-Wasserstein distance on the real line (and in particular, on the unit interval) is the $L_1$-distance of the cumulative distribution functions, that is, 
\begin{align}
d_{\cW_1([0,1],|.|)}(\mu, \nu)=\int_{[0,1]} |F_\mu(x)-F_\nu(x)| ~\dd x.
    \nonumber
\end{align}
Consequently, any flow of CDFs interpolating between $F_\mu$ and $F_\nu$ gives rise to a $\cW_1$-geodesic between $\mu$ and $\nu,$ and the interpolation is never unique, unless $\nu=\mu.$ By $\cW_1(Y,r) \simeq \cW_1([0,1], |.|),$ this implies that $\cW_1(Y,r)$ also has the property that any pair of distinct points in it are connected by multiple geodesics. 
Consequently, for any $a,b \in Y, a \neq b$ there is a point "strictly between them", that is, there exists a $z \in Y \setminus \{a,b\}$ such that $r(a,z)+r(z,b)=r(a,b).$ Indeed, assume by contradiction that $r(a,z)+r(z,b)>r(a,b)$ for all $z \notin \{a,b\}.$ For any probability measure $\mu,$ 
\begin{align}
    d_{\cW_1(Y,r)}(\delta_a, \mu)+d_{\cW_1(Y,r)}(\mu, \delta_b)
    =\int_{Y} (r(a,z)+r(z,b)) ~\dd \mu(z),
    \nonumber
\end{align}
which is equal to $d_{\cW_1(Y,r)}(\delta_a, \delta_b)=r(a,b)$ if and only if $\mu$ is supported on $\{a,b\}.$ Consequently, the standard convex combination $[0,1] \ni \lambda \mapsto (1-\lambda)\delta_a+\lambda \delta_b$ is the unique geodesic connecting $\delta_a$ and $\delta_b,$ a contradiction. 
\par
And this forces $(Y,r)$ to be geodesic.
Indeed, according to \cite[Definition 1]{Boyd-Wong-PAMS-69} (see also \cite[p. 41]{Blumenthal-distance-geometry-book}), the metric space $(Y,r)$ is \emph{metrically convex}. It is also complete because it is compact, and hence, by \cite[Lemma 1]{Boyd-Wong-PAMS-69} (see also \cite[Theorem 14.1]{Blumenthal-distance-geometry-book}), this implies that for any $\alpha \in (0,1)$ and any $a,b \in (Y,r)$ there exists a $z=z(\alpha)$ such that 
\begin{align}
    r(a,z)=(1-\alpha)r(a,b) \text{ and } r(z,b)=\alpha r(a,b).
\end{align}
In particular, midpoints exist. Iterating midpoint insertion constructs an isometric para\-metrization of the dyadic rational points of $[0,r(a,b)].$ By the completeness of $(Y,d_Y),$ we can extend this parametrization to the whole interval, which gives the required geodesics. 
\par
Since $(Y,r)$ is compact, geodesic, and has diameter $1,$ it contains an isometric copy of $([0,1],|.|).$ Let $\gamma: ([0,1],|.|) \to (Y,r)$ denote this isometric embedding. By push-forward, this induces an isometric embedding $\gamma_{\#}$ of $\cW_1([0,1],|.|)$ into $\cW_1(Y,r).$ Composing this embedding with the isometry witnessing \eqref{eq:W-1-Y-r-isom-W-1-01-std}, we get an isometric embedding of $\cW_1([0,1],|.|)$ into itself. But isometric self-embeddings of a compact metric space are surjective, and hence so is $\gamma_{\#}.$ But all measures in the image of $\gamma_{\#}$ are supported on $\gamma([0,1]).$ The surjectivity of $\gamma_{\#},$ applied to $\delta_y$ for each $y \in Y$ forces $Y=\gamma(I),$ which means that $(Y,r) \simeq ([0,1],|.|),$ and taking the $p$-th roots yields the desired conclusion 
\begin{align}
    %\label{eq:des-conc}
    (Y,d_Y) \simeq ([0,1],|.|^{1/p}). \nonumber
\end{align}
\par
Theorem \ref{thm:recovery-is-generic} states that the Dirac preservation property is generic in $\mathfrak{M},$ and hence so is the recovery property.

The proof strategy for Theorem \ref{thm:recovery-is-generic} is similar to that of Theorem \ref{thm:rigidity-is-generic}. We show that the 
set of metric spaces satisfying the Dirac preservation property is a dense $G_{\delta}$ set in $\mathfrak{M}$. The proof again has two parts. The density property will follow by looking at finite anisometric metric spaces, while the $G_{\delta}$ property will follow from the 
continuity of the map 
$$T_p : \mathfrak{M} \to \mathfrak{M}, \ X \ \to \ \Wp(X,d_X).$$

The following statement shows the first step in this scheme: the density of spaces satisfying the Dirac preservation and hence the recovery property. 

\begin{proposition} \label{prop:isom-underlying-spaces-finite} 
Let $(X,d_X)$ be a finite, totally anisometric metric space, let $p>1,$ and assume that $\Wp(X,d_X)$ is isometric to $\Wp(Y,d_Y)$ for some Polish metric space $(Y, d_Y).$ Then for any isometry $\Phi: \Wp(X,d_X) \to \Wp(Y,d_Y),$ the image $\Phi(\delta_x) \in \Wp(Y,d_Y)$ is a Dirac measure for all $x \in X,$ and for every $y \in Y$ there is an $X \in X$ such that $\Phi(\delta_x)=\delta_y.$ Consequently, $Y$ is isometric to $X.$ 
\end{proposition} 

\begin{proof}
Let $n:=|X|,$ and let $\Phi: \Wp(X,d_X) \to \Wp(Y,d_Y)$ be an isometry. The topological dimension of $\cW_p(X,d_X)$ is $n-1,$ and hence so is the dimension of $\cW_p(Y,d_Y),$ which implies that $|Y|=n$ as well. 
First we show that $\Phi(\delta_x)$ is a Dirac mass for all $x\in X.$ To do so, sssume that the unique pair of points realizing the diameter in $(X,d_X)$ is $\{x_1, x_n\}.$ Then the unique pair of point realizing the diameter in $\cW_p(X,d_X)$ is $\{\delta_{x_1},\delta_{x_n}\}.$ Consequently, $\{\Phi(\delta_{x_1}),\Phi(\delta_{x_n})\}$ is the unique pair of point realizing the diameter in $\cW_p(Y,d_Y),$ which implies that both $\Phi(\delta_{x_1})$ and $\Phi(\delta_{x_n})$ are Dirac masses. Let $y_1, y_n \in Y$ be defined by $\delta_{y_1}=\Phi(\delta_{x_1})$ and $\delta_{y_n}=\Phi(\delta_{x_n}).$ Let $\sim$ denote the binary relation on $\cW_p(X,d_X)$ and $\cW_p(Y,d_Y)$ characterized by Lemma \ref{lemma:metric-ratio-singleton}. Let us introduce the sets $\mathcal{M}:=\lers{\mu \in \cW_p(X,d_X) \, \middle| \, \mu \sim \delta_{x_1}}$ and $\mathcal{N}:=\lers{\nu \in \cW_p(Y,d_Y) \, \middle| \, \nu \sim \delta_{y_1}}.$ By point (a) of Lemma \ref{lemma:metric-ratio-singleton}, 
\begin{align} \label{eq:M-char-2}
    \mathcal{M}=\lers{(1-c)\delta_{x_1} + c \delta_x \, \middle| \, (x_1 ,x) \in \cF, \ c \in [0,1]},
\end{align}
and similarly, 
\begin{align} \label{eq:N-char}
    \mathcal{N}=\lers{(1-c)\delta_{y_1} + c \delta_y \, \middle| \, (y_1 ,y) \in \cF, \ c \in [0,1]}.
\end{align}
The relation $\sim$ is characterized in point (b) of Lemma \ref{lemma:metric-ratio-singleton} purely in terms of the $p$-Wasserstein metric, and hence $\mathcal{N}=\{ \Phi(\mu) \, | \, \mu \in \mathcal{M}\}.$ Clearly, $d_{\Wp}(\delta_{x_1}, (1-c)\delta_{x_1} + c \delta_x)=c^{1/p} d_X(x_1,x),$ and hence $\max \{d_{\Wp}(\delta_{x_1}, \mu) \, | \, \mu \in \mathcal{M} \}=\max \{d_X(x_1,x) \,|\, (x_1,x) \in \cF\},$ and 
\begin{align}
    \argmax \{d_{\Wp}(\delta_{x_1}, \mu) \, | \, \mu \in \mathcal{M} \}=\delta_{x_2},
\end{align}
where $x_2$ is defined by $x_2:=\argmax \{ d_X(x_1,x) \,|\, (x_1,x) \in \cF \}. $ Similarly, 
\begin{align}
    \argmax \{d_{\Wp}(\delta_{y_1}, \nu) \, | \, \mu \in \mathcal{N} \}=\delta_{y_2},
\end{align}
where $y_2$ is defined by $y_2:=\argmax \{ d_Y(y_1,y) \,|\, (y_1,y) \in \cF \}.$ The isometry $\Phi$ preserves the unique maximal distance, and hence $\Phi(\delta_{x_2})=\delta_{y_2}.$ Now let 
\begin{align}
    \mathcal{M}_1:=\mathcal{M} \setminus \lers{(1-c)\delta_{x_1} + c \delta_{x_2} \, \middle| \, c \in (0,1]}.
\end{align}
The metric characterization of $\sim,$ that is, point (b) in Lemma \ref{lemma:metric-ratio-singleton} together with $\Phi(\delta_{x_1})=\delta_{y_1}$ and $\Phi(\delta_{x_2})=\delta_{y_2}$ implies that $\Phi\ler{(1-c)\delta_{x_1}+c \delta_{x_2}}=(1-c)\delta_{y_1}+c \delta_{y_2}$ for all $c \in [0,1].$ Therefore, for the set 
\begin{align}
    \mathcal{N}_1:=\mathcal{N} \setminus \lers{(1-c)\delta_{y_1} + c \delta_{y_2} \, \middle| \, c \in (0,1]}
\end{align}
we have that 
\begin{align}
    \lers{\Phi(\mu) \, | \, \mu \in \mathcal{M}_1}=\mathcal{N}_1. 
\end{align}
Consequently, with the notation $x_3:=\argmax \{ d_X(x_1,x) \,|\, (x_1,x) \in \cF, x \neq x_2\}$ we have 
\begin{align}
    \argmax \{d_{\Wp}(\delta_{x_1}, \mu) \, | \, \mu \in \mathcal{M}_1 \}=\delta_{x_3},
\end{align} 
and similarly, with 
$y_3:=\argmax \{ d_Y(y_1,y) \,|\, (y_1,y) \in \cF, y \neq y_2\}$ we have 
\begin{align}
    \argmax \{d_{\Wp}(\delta_{y_1}, \nu) \, | \, \nu \in \mathcal{N}_1 \}=\delta_{y_3},
\end{align} 
which implies that $\Phi(\delta_{x_3})=\delta_{y_3}.$ And so on, one can conclude that for all $x \in X$ with $(x,x_1) \in \cF$ there is a $y \in Y$ such that $\Phi(\delta_{x})=\delta_y.$ 
Now consider the $\cF$-neighbors of $x_2$ that have not been covered so far, that is, the set 
\begin{align}
   \lers{x\in X \, \middle| \, (x_2,x) \in \cF \text{ and } (x_1, x) \notin \cF },
\end{align}
and the corresponding set of measures 
\begin{align}
     \mathcal{R}:=\lers{(1-c) \delta_{x_2}+c \delta_x \, \middle| \, (x_2,x) \in \cF \text{ and } (x_1, x) \notin \cF, \, c \in [0,1]}.
\end{align}
Playing the same game with $\mathcal{R}$ that we played with $\mathcal{M}$ we can prove for newer points in $X$ that the Dirac mass supported on them is sent to a Dirac mass by $\Phi.$ 
The graph of the relation $\cF$ is connected on $X,$ see Lemma \ref{lemma:X-F-graph-connected}, and hence we conclude that for every $x \in X$ there is a $y \in Y$ such that $\Phi(\delta_x)=\delta_y.$ This means that there is a map $F: X \to Y$ such that $\Phi(\delta_x)=\delta_{F(x)}$ for all $x \in X,$ and the usual 
\begin{align}
    d_Y(F(x),F(x'))=d_{\cW_p}(\delta_{F(x)},\delta_{F(x')})
=d_{\cW_p}(\Phi(\delta_x),\Phi(\delta_{x'}))
=d_{\cW_p}(\delta_x,\delta_{x'})
=d_X(x,x') \nonumber
\end{align}
argument demonstrates that $F$ is an isometric embedding of $X$ into $Y.$ It follows from $|X|=|Y|=n$ that $F$ is also surjective, which implies that for all $y \in Y$ there is an $x \in X$ such that $\Phi(\delta_x)=\delta_y,$ and that $Y$ is isometric to $X.$
\end{proof}

The next proposition shows that the set of compact metric spaces with the Dirac preservation property form a $G_{\delta}$ set in $\mathfrak{M}$. The combination of Proposition \ref{prop:isom-underlying-spaces-finite} and Proposition \ref{prop:G-delta-rec} makes the proof of Theorem \ref{thm:recovery-is-generic} complete.

\begin{proposition} \label{prop:G-delta-rec}
For $p\geq 1$ and $\varepsilon >0$ the set 
\begin{align*} %\label{def: Ueps}
U_{\varepsilon}:= \{ X \in \mathfrak{M}: \ \text{if} \ Y\in \mathfrak{M }\ \text{and}  
\ \Phi\in \ \mathrm{Isom}(\Wp(X,d_X), \Wp(Y,d_Y)) \ \Rightarrow \forall x \in X \ \exists \ y \in Y \\
~\text{with}~ \ d_{\Wp}(\Phi(\delta_x), \delta_y) < \varepsilon \ ~\text{and}~ \  \forall \ y \in Y , \exists \ x \in X ~\text{with}~ \ d_{\Wp}(\Phi(\delta_x), \delta_y) < \varepsilon\}
\end{align*}
is open in $\mathfrak{M}$. 
Consequently the set 
\begin{align*} 
 \mathcal{R}= \{ X \in \mathfrak{M}:  \text{if} \ Y\in \mathfrak{M }\ ~\text{and}~  
\ \Phi\in \ \mathrm{Isom}(\Wp(X,d_X), \Wp(Y,d_Y)) \ \Rightarrow \forall x \in X \ \exists \ y \in Y \\
~\text{with}~ \ \Phi(\delta_x)= \delta_y \ ~\text{and}~ \  \forall \ y \in Y , \exists \ x \in X ~\text{with}~ \ \Phi(\delta_x)= \delta_y\} = \bigcap_{n=1}^{\infty} U_{\frac{1}{n}}
\end{align*}
is a $G_{\delta}$ set in $\mathfrak{M}$.
\end{proposition}

\begin{proof} The idea of the proof is similar to the one of Lemma \ref{lemma:R-G-delta}. We shall show that the complement $\mathcal{Q}_{\varepsilon} = \mathfrak{M} \setminus U_{\varepsilon}$ is closed in $\mathfrak{M}$. 
Let us note that $\mathcal{Q}_{\varepsilon}$ can be written as 
\begin{align*} %\label{def:Q_eps}
\mathcal{Q}_{\varepsilon}:= \{ X \in \mathfrak{M}: ~&\exists  \ Y\in \mathfrak{M }\ \text{and}  
\ \Phi\in \ \mathrm{Isom}(\Wp(X,d_X), \Wp(Y,d_Y)) \\
\text{such that} \ &\exists x \in X \ \text{such that} \ \forall \ y \in Y 
 \ d_{\Wp}(\Phi(\delta_x), \delta_y) \geq \varepsilon \\
 \text{or} \  &\exists  y \in Y , \ \text{such that} \  \forall \ x \in X  \ d_{\Wp}(\Phi(\delta_x), \delta_y) \geq \varepsilon\}.
\end{align*}

To show that $\mathcal{Q}_{\varepsilon}$ is closed let us consider a sequence $(X_k)_{k\in\mathbb{N}} \subset \mathcal{Q}_{\varepsilon}$  converging to some $X \in \mathfrak{M}$ in the Gromov-Hausdorff sense. We have to show that $X \in \mathcal{Q}_{\varepsilon}$. 

Without loss of generality we assume that that all elements of the sequence as well as $X$ are contained in the same space $Z \in \mathfrak{M}$.

Let us assume that for infinitely many indices $k$ the first option holds:  there exist $Y_k$ 
and $\Phi_k: \Wp(X_k,d_{X_k}) \to \Wp(Y_k,d_{Y_k})$ and $ x_k \in X_k $ such that for all  $ y \in Y_k $
we have that  
$$ 
d_{\Wp}(\Phi_k(\delta_{x_k}), \delta_y) \geq \varepsilon.
$$ 
By the continuity of the Wasserstein map \eqref{eq:1-lipshitz}, 
\begin{align}
\label{eq:Wass-map-limit}
\lim_{k \to \infty } \cW_p(Y_k,d_{Y_k})=\lim_{k \to \infty } \cW_p(X_k,d_{X_k})=\cW_p(X,d_X).
\end{align}
In particular, $(\cW_p(Y_k,d_{Y_k}))_{k\in\mathbb{N}}$ converges to a compact space, hence the set $\lers{\cW_p(Y_k,d_{Y_k})}$ is uniformly totally bounded (see \cite[Def. 7.4.13.]{BBI}) meaning that its elements have uniformly bounded diameters and uniform covering bounds at every positive scale. As $\cW_p(Y_k,d_{Y_k})$ contains an isometric copy of $Y_k,$ the sequence of base spaces $(Y_k)_{k\in\mathbb{N}}$ inherits these uniform bounds. Therefore, the set $\{Y_k\}_{k \in \N}$ is uniformly totally bounded, and hence by \cite[Thm. 7.4.15.]{BBI} pre-compact, which means that  $(Y_k)_{k\in\mathbb{N}}$ has a convergent subsequence converging to a limit $Y \in \mathfrak{M}$ in the Gromov-Hausdorff distance. To simplify notation we denote this further subsequence by $(X_k)_{k\in\mathbb{N}}$, $(Y_k)_{k\in\mathbb{N}}$ and isometry $\Phi_k: \Wp(X_k,d_{X_k}) \to \Wp(Y_k,d_{Y_k})$ with the property that there exists $x_k \in X_k$ such that 

$$d_{\Wp}(\Phi_k(\delta_{x_k}), \delta_y) \geq \varepsilon $$
for all $y\in Y_k$. Note that $X_k \to X$, $Y_k \to Y$ in the Gromov-Hausdorff distance as $k\to \infty$. Applying  \eqref{eq:1-lipshitz} we obtain that $\Wp(X_k,d_{X_k}) \to \Wp(X,d_X) $ and $ \Wp(Y_k,d_{Y_k}) \to \Wp(Y,d_Y)$ in $\mathfrak{M}$ as $k\to \infty$. Using the same idea as in Lemma \ref{lemma:R-G-delta} we can select a further subsequence (denoted again in the same way) such that $\Phi_k \to \Phi$, where $\Phi: X \to Y$ is an isometry. Selecting one more time a further subsequence we can assume without loss of generality that $x_k \to x$ where $x\in X$.

Passing to the limit in the inequality $ d_{\Wp}(\Phi_k(\delta_{x_k}), \delta_y) \geq \varepsilon$
as $k\to \infty$ we obtain that $$ d_{\Wp}(\Phi(\delta_{x}), \delta_y) \geq \varepsilon$$ for all $y \in Y$, showing that 
the limit space $X \in \mathcal{Q}_{\varepsilon}$. 

Assume now the second option: for infinitely many indices $k$  there exist $Y_k\in \mathfrak{M}$ and $y_k \in Y_k$
and $\Phi_k: \Wp(X_k,d_{X_k}) \to \Wp(Y_k,d_{Y_k})$ such that for all  $ x \in X_k $
we have that  
$$ d_{\Wp}(\Phi_k(\delta_{x}), \delta_{y_k}) \geq \varepsilon .$$

We have to show again, that the limit space $X$ is in $\mathcal{Q}_{\varepsilon}$ This is done again by producing a limiting space $Y$ and isometry $\Phi: \Wp(X,d_X) \to \Wp(Y,d_Y)$ and a point $y\in Y$ such that $d_{\Wp}(\Phi(\delta_{x}), \delta_{y}) \geq \varepsilon$ for all $x \in X$. The proof of this is analogous to the one from the previous case. We leave the details to the interested reader. 
\end{proof}

\section*{AI Disclosure Statement.} We used OpenAI's GPT-6 Astra large language model to draft the proof of Proposition \ref{prop:rigid-G-delta} and Lemma \ref{lemma:R-G-delta}, and to find the argument showing that the property DPPp is strictly stronger than RECp. The authors take full responsibility for the results presented in this paper.

\end{document}